\documentclass[final,1p,times]{elsarticle}

\usepackage{amssymb}
 \usepackage{amsthm}
\usepackage{amscd}
\usepackage{bm}
\usepackage{color}
\usepackage{amsmath}
\usepackage{amsfonts}
\usepackage{amssymb}
\usepackage{graphicx}
\newtheorem{theorem}{Theorem}

\newtheorem{lemma}[theorem]{Lemma}

\usepackage{mathrsfs}
\usepackage{titletoc}
\usepackage{bbding}

\newcommand{\ra}{\rightarrow}
\newcommand{\p}{\partial}
\newcommand{\f}{\frac}

\newcommand{\be}{\begin{equation}}
\renewcommand{\ra}{\rightarrow}
\newcommand{\ee}{\end{equation}}
\newcommand{\bea}{\begin{eqnarray}}
\newcommand{\eea}{\end{eqnarray}}
\newcommand{\bna}{\begin{eqnarray*}}
\newcommand{\ena}{\end{eqnarray*}}

\renewcommand{\le}{\left}
\newcommand{\ri}{\right}

\journal{***}

\begin{document}

\begin{frontmatter}

\title{Mean field type equations on line bundle\\ over a closed Riemann surface \tnoteref{sw}}

\author[ruc]{Jie Yang}
\ead{yangjie1112@ruc.edu.cn}

\author[ruc]{Yunyan Yang\corref{jy}}
\address[ruc]{Department of Mathematics, Renmin University of China, Beijing 100872, China}
\ead{yunyanyang@ruc.edu.cn}

\cortext[jy]{Corresponding author.}
\tnotetext[sw]{This research is partly supported by the National Natural Science Foundation of China (Grant No. 11721101),
and by the National Key Research and Development Project SQ2020YFA070080.}

\begin{abstract}

Let $(\mathcal{L},\mathfrak{g})$ be a line bundle over  a closed Riemann surface $(\Sigma,g)$,
$\Gamma(\mathcal{L})$ be the set of all smooth sections,
and $\mathcal{D}:\Gamma(\mathcal{L})\ra T^\ast\Sigma\otimes \Gamma(\mathcal{L})$ be a connection independent of the bundle
metric $\mathfrak{g}$, where $T^\ast\Sigma$ is the cotangent bundle.
Suppose that there exists a global unit frame $\zeta$ on $\Gamma({\mathcal{L}})$. Precisely for
any $\sigma\in\Gamma(\mathcal{L})$, there exists a unique smooth function $u:\Sigma\ra\mathbb{R}$ such that
$\sigma=u\zeta$ with $|\zeta|\equiv 1$ on $\Sigma$.
 For any real number $\rho$, we define a functional
$\mathcal{J}_\rho:W^{1,2}(\Sigma,\mathcal{L})\ra\mathbb{R}$ by
$$\mathcal{J}_\rho(\sigma)=\f{1}{2}\int_\Sigma|\mathcal{D} \sigma|^2dv_g+\f{\rho}{|\Sigma|}\int_\Sigma\langle\sigma,\zeta\rangle dv_g
-\rho\log\int_\Sigma h e^{\langle\sigma,\zeta\rangle}dv_g,$$
where $W^{1,2}(\Sigma,\mathcal{L})$ is a completion of
$\Gamma(\mathcal{L})$ under the usual Sobolev norm, $|\Sigma|$ is the area of $(\Sigma,g)$, $h:\Sigma\ra\mathbb{R}$ is a strictly positive smooth function
and  $\langle\cdot,\cdot\rangle$ is the inner product induced by $\mathfrak{g}$.
The Euler-Lagrange equations of $\mathcal{J}_\rho$ are called mean field type equations.
Write $\mathcal{H}_0=\{\sigma\in W^{1,2}(\Sigma,\mathcal{L}):\mathcal{D}\sigma=0\}$ and
$$\mathcal{H}_1=\le\{\sigma\in W^{1,2}(\Sigma,\mathcal{L}):\int_\Sigma \langle\sigma,\tau\rangle dv_g=0,\,\,\forall \tau \in \mathcal{H}_0\ri\}.$$
Based on the variational method, we prove that
$\mathcal{J}_\rho$ has a constraint critical point on the space $\mathcal{H}_1$ for any $\rho<8\pi$;
Based on blow-up analysis, we calculate the exact value of $\inf_{\sigma\in\mathcal{H}_1}\mathcal{J}_{8\pi}(\sigma)$,
provided that it is not achieved by any $\sigma\in\mathcal{H}_1$;  If we further assume $\mathcal{D}\zeta=0$,
$\mathcal{J}_\rho:W^{1,2}(\Sigma,\mathcal{L})\ra\mathbb{R}$ is reduced to a functional related to
the classical mean field equation.
\end{abstract}


\end{frontmatter}

\section{Introduction}
Let $(\Sigma , g )$ be a closed Riemann surface. The well-known mean field equation has aroused the interests of
many mathematicians for a long time. It is written as
\begin{eqnarray}\label{5}
\Delta_g u=\rho\left(\frac{h e^u}{\int_{\Sigma} h e^u\mathrm{d}v_g}-\frac{1}{|\Sigma|}\right),
\end{eqnarray}
where $\Delta_g$ is the Laplace-Beltrami operator,  $\rho$ is a real number, $h:\Sigma\ra\mathbb{R}$ is a function,
and $|\Sigma|$ denotes the area of $\Sigma$.
This equation appears in the prescribed Gaussian curvature problem in conformal geometry
 \cite{KW,Chang-Yang1,Chang-Yang2,S}
and also in the abelian Chern-Simons-Higgs model in physics \cite{Caffarelli,DJLW98,DJLW99B,Tarantello, YangYisong}.
Let $W^{1,2}(\Sigma)$ be the standard Sobolev space with respect to the norm
\begin{eqnarray*}\label{a0}
\|u\|_{W^{1,2}(\Sigma)}=\left(\int_\Sigma(|\nabla_g u|^2+|u|^2)dv_g\right)^{1/2},
\end{eqnarray*}	
where $\nabla_g$ stands for the gradient operator on $(\Sigma , g )$.
Note that \eqref{5} has a variational structure. In particular, solutions of \eqref{5} are critical points of the functional
$J_{\rho}:W^{1,2}(\Sigma)\ra\mathbb{R}$, which is defined by
\begin{eqnarray}\label{Jrho}
J_{\rho}(u)=\frac{1}{2}\int_{\Sigma}|\nabla_g u|^{2} \mathrm{d}v_g+\f{\rho}{|\Sigma|}\int_\Sigma{u}dv_g -\rho \log\int_{ \Sigma} h e^{u} \mathrm{d}v_g.
\end{eqnarray}

The solvability of the mean field equation (\ref{5}) is closely related to the Trudinger-Moser inequality
contributed by Moser \cite{Moser71}, Aubin \cite{Aubin}, Fontana \cite{Fontana}, among others.
Such an inequality together with a direct method of variation implies that  $J_{\rho}$  has a critical point when $\rho < 8\pi$.
In the critical case $\rho= 8\pi$, under certain geometric and analytic assumptions,
a minimizing solution was found by Ding-Jost-Li-Wang \cite{DJLW} via the method of blow-up analysis. This result
 was extended to various situations, see for examples \cite{Li-Zhu,Wang,NT,Yang-Zhu-Proc,Yang-Zhu-Sci, Sun-Zhu,Sun-Z}. For the related Nirenberg problem, the readers are referred to \cite{Aubin,Chang-Yang1,Chang-Yang2,KW,Moser71, Moser73,S}
and the references therein. However, the case $\rho>8\pi$ is much different from $\rho\leq 8\pi$.
Using a min-max method, Ding-Jost-Li-Wang \cite{DJLW99} obtained that $J_{\rho}$ has a critical point if $\rho\in (8 \pi, 16 \pi)$, $h>0$ and the genus of the surface is nonzero. A topological method was introduced by Djadli-Malchiodi \cite{DM} and Djadli \cite{Djadli} to
solve this problem. In particular, for any $\rho \in(8k \pi, 8(k+1) \pi)$, $ k \in \mathbb{N}^{*}$,  $J_\rho$ has a critical point with no
assumption on the topology of the surface. Along this direction, there had been a lot of work such as  \cite{AJY,Malchiodi2015,MN,Malchiodi-DCDS, MR}
and the references therein. This problem still arouses many people's interest \cite{JDE20,changesign,Sun-Wang-Yang,Li-Sun-Yang}.

In the current paper, our aim is to extend the mean field equation from Riemann surface to the line bundle setting.
Given a line bundle $(\mathcal{L},\mathfrak{g})$ over a closed Riemann surface $(\Sigma,g)$, $\mathfrak{g}$ is a bundle metric.
Let $\Gamma(\mathcal{L})$ be the set of all smooth sections $\sigma:\Sigma\ra\mathcal{L}$, $T^\ast\Sigma$ be the cotangent bundle, and
$\mathcal{D}:\Gamma(\mathcal{L})\ra T^\ast\Sigma\otimes \Gamma(\mathcal{L})$ be a connection. Here and in the sequel we do not assume $\mathfrak{g}$ is
parallel with respect to the connection $\mathcal{D}$.
 Define a Sobolev space $W^{1,2}(\Sigma,\mathcal{L})$ as a completion of $\Gamma(\mathcal{L})$ under the norm
$$\|\sigma\|_{W^{1,2}(\Sigma,\mathcal{L})}=\le(\int_\Sigma(|\mathcal{D}\sigma|^2+|\sigma|^2)dv_g\ri)^{1/2},$$
where both of the metrics on $\Gamma(\mathcal{L})$ and $T^\ast\Sigma\otimes\Gamma(\mathcal{L})$ are denoted by the same symbol $|\cdot|$.
A unit section $\zeta\in \Gamma(\mathcal{L})$ is said to be a global unit frame provided that $|\zeta(x)|\equiv 1$ for all $x\in\Sigma$ and for any fixed $\sigma\in \Gamma(\mathcal{L})$, there exists a unique
smooth function
$u:\Sigma\ra\mathbb{R}$ such that $\sigma(x)=u(x)\zeta(x)$ for all $x\in\Sigma$. It is easy to see that any global unit frame
 on $\Gamma(\mathcal{L})$ would be also a global unit frame on $W^{1,2}(\Sigma,\mathcal{L})$.
 Then we naturally extend the connection $\mathcal{D}$ from $\Gamma(\mathcal{L})$ to $W^{1,2}(\Sigma,\mathcal{L})$.
Indeed, there exists a unique $\omega\in \Gamma(T^\ast\Sigma)$ satisfying $\mathcal{D}\zeta=\omega\otimes\zeta$.
As a consequence, for any $\sigma=u\zeta\in W^{1,2}(\Sigma,\mathcal{L})$, we may write by the Leibniz rule
\be\label{Lebn}\mathcal{D}\sigma=du\otimes\zeta+u\omega\otimes\zeta.\ee
 Then for any real number $\rho$, with the help of such a $\zeta$, we may define a functional
$\mathcal{J}_\rho:W^{1,2}(\Sigma,\mathcal{L})\ra\mathbb{R}$ by
\be\label{J-funct}\mathcal{J}_\rho(\sigma)=\f{1}{2}\int_\Sigma|\mathcal{D} \sigma|^2dv_g+\f{\rho}{|\Sigma|}\int_\Sigma\langle\sigma,\zeta\rangle dv_g-\rho\log\int_\Sigma he^{\langle\sigma,\zeta\rangle}dv_g,\ee
where $\langle\cdot,\cdot\rangle$ denotes the inner product induced by $\mathfrak{g}$, and $h:\Sigma\ra\mathbb{R}$ is a positive smooth function.
It should be noticed that the above functional $\mathcal{J}_\rho$ depends on the global unit frame $\zeta$.
Define two section sets
\be\label{H0}\mathcal{H}_0=\le\{\sigma\in W^{1,2}(\Sigma,\mathcal{L}):\mathcal{D}\sigma=0\ri\}\ee
and
\be\label{H1}\mathcal{H}_1=\le\{\sigma\in W^{1,2}(\Sigma,\mathcal{L}):\int_\Sigma \langle\sigma,\tau\rangle dv_g=0,\,\,\forall \tau \in \mathcal{H}_0\ri\}.\ee
Here and throughout this paper we abuse some notations slightly. Both zero section and real number zero are denoted by the same $0$. However, readers can recognize them from the context.
Obviously $\mathcal{H}_0\cap \mathcal{H}_1=\{0\}$ and
\bna W^{1,2}(\Sigma,\mathcal{L})=\mathcal{H}_0\oplus \mathcal{H}_1.\ena

In the line bundle setting, we restrict
sections on the space $\mathcal{H}_1$ mainly because the Poincar\'e inequality holds there. To write explicitly the Euler-Lagrange equations for
minimizers of $\mathcal{J}_\rho$, we also need to understand $\mathcal{H}_0$. We will prove later $\mathcal{H}_0$ is a finite
dimensional linear space, in particular ${\rm dim}\mathcal{H}_0\leq 1$. As for the dimension of $\mathcal{H}_0$, there are two possibilities: ${\rm dim} \mathcal{H}_0=0$
or ${\rm dim} \mathcal{H}_0= 1$. In view of (\ref{Lebn}) and the definitions of $\mathcal{D}$ and $\zeta$, for any $\tau=u\zeta\in\mathcal{H}_0$, we have
$du\otimes\zeta+u\omega\otimes\zeta=(du+u\omega)\otimes\zeta=0$,
or equivalently \be\label{zero1}du+u\omega= 0\quad{\rm on}\quad\Sigma.\ee
Obviously $u\equiv 0$ is a solution of (\ref{zero1}). If (\ref{zero1}) has no other solution, then
$\mathcal{H}_0=\{0\}$, and thus ${\rm dim}\mathcal{H}_0=0$; If (\ref{zero1}) has a solution $u\not\equiv 0$, then
${\rm dim}\mathcal{H}_0= 1$. Let us give an example. Assume $(\Sigma,g)$ is a closed Riemann surface, $\mathcal{L}=\Sigma\times\mathbb{R}$ is a
line bundle
and $\mathfrak{g}$ is the standard Euclidean metric on $\mathbb{R}$.
Take a global unit frame $\zeta(x)=(x,1)$ for all $x\in\Sigma$ and a bundle connection $\mathcal{D}$ such that
$\mathcal{D}\zeta=\omega\otimes\zeta=0$, which implies $\omega=0$. Obviously $W^{1,2}(\Sigma,\mathcal{L})\cong W^{1,2}(\Sigma)$,
and (\ref{zero1}) has only solutions $u\equiv C$ for all $C\in\mathbb{R}$. Hence $\mathcal{H}_0\cong \mathbb{R}$ and
${\rm dim}\mathcal{H}_0=1$.
In the case ${\rm dim}\mathcal{H}_0= 1$, we may assume
$\mathcal{H}_0={\rm span}\{\tau_1\}$,
where $\{\tau_1\}$ is a normal basis, in particular
$$\int_\Sigma \langle\tau_1,\tau_1\rangle dv_g=1.
$$

If $\sigma$ is a minimizer of $\mathcal{J}_\rho$ on $\mathcal{H}_1$, then it satisfies
\be\label{H0-case}\Delta_{\mathcal{L}}\sigma=\rho\f{ he^{\langle\sigma,\zeta\rangle}\zeta}{\int_\Sigma he^{\langle\sigma,\zeta\rangle}dv_g}-
\f{\rho}{|\Sigma|}\zeta\ee
in the case ${\rm dim}\mathcal{H}_0=0$, and
\be\label{E-L}\le\{\begin{array}{lll}\Delta_{\mathcal{L}}\sigma=\rho\f{ he^{\langle\sigma,\zeta\rangle}\zeta}{\int_\Sigma he^{\langle\sigma,\zeta\rangle}dv_g}-
\f{\rho}{|\Sigma|}\zeta-\lambda_1\tau_1\\[1.5ex]
\lambda_1=\rho\f{\int_\Sigma he^{\langle\sigma,\zeta\rangle}\langle\zeta,\tau_1\rangle dv_g}{\int_\Sigma he^{\langle\sigma,\zeta\rangle} dv_g}
-\f{\rho}{|\Sigma|}\int_\Sigma \langle \zeta, \tau_1\rangle dv_g
\end{array}\ri.\ee
in the case ${\rm dim}\mathcal{H}_0=1$,
where $\Delta_{\mathcal{L}}=\mathcal{D}^\ast\mathcal{D}$ is the bundle Laplace-Beltrami operator,
and $\mathcal{D}^\ast$ is the dual operator of $\mathcal{D}$, which is defined by
\bna\label{Dstar}\int_\Sigma \langle\mathcal{D}^\ast\mathcal{D}\sigma,\phi\rangle dv_g=\int_\Sigma \langle\mathcal{D}\sigma,\mathcal{D}\phi\rangle dv_g,
\quad\forall \phi\in W^{1,2}(\Sigma,\mathcal{L}).\ena
Let us explain the geometric meaning of the equation (\ref{E-L}). Since $\mathcal{H}_i$ is a subset of $L^2(\Sigma,\mathcal{L})$,
and $W^{1,2}(\Sigma,\mathcal{L})$ is dense in $L^2(\Sigma,\mathcal{L})$ under the usual $L^2$ norm, we have $L^2(\Sigma,\mathcal{L})=\overline{\mathcal{H}}_0\oplus\overline{\mathcal{H}}_1$, where $\overline{\mathcal{H}}_i$ is a closure of $\mathcal{H}_i$
in $L^2(\Sigma,\mathcal{L})$, $i=1,2$. This is an orthogonal decomposition of $L^2(\Sigma,\mathcal{L})$.
There exists a natural projection $\mathcal{P}: L^2(\Sigma,\mathcal{L})\ra\overline{\mathcal{H}}_1$, which maps a $\xi\in L^2(\Sigma,\mathcal{L})$
to
$$\mathcal{P}(\xi)=\xi-\le(\int_\Sigma \langle\xi,\tau_1\rangle dv_g\ri)\tau_1.$$
In this framework, the equation (\ref{E-L}) can be written as
$$\Delta_{\mathcal{L}}\sigma=\mathcal{P}\le\{\rho\le(\f{he^{\langle\sigma,\zeta\rangle}}{\int_\Sigma he^{\langle\sigma,\zeta\rangle}dv_g}
-\f{1}{|\Sigma|}\ri)\zeta\ri\}.$$

Concerning critical points of $\mathcal{J}_\rho$, we  have the following:
\begin{theorem}\label{T1}
Let $(\mathcal{L},\mathfrak{g})$ be a line bundle over a closed Riemann surface $(\Sigma,g)$, $\mathcal{D}$ be its bundle connection.
Suppose there exists a global unit frame $\zeta$ on $\Gamma(\mathcal{L})$.
Let $\mathcal{J}_\rho$ and $\mathcal{H}_1$ be defined as in (\ref{J-funct}) and (\ref{H1}) respectively.
 Then for any $\rho<8\pi$, $\mathcal{J}_\rho$ has a minimizer $\sigma_\rho$ on
$\mathcal{H}_1$.
\end{theorem}

The proof of Theorem \ref{T1} is based on a method of variation. It should be remarked that in this theorem,
the infimum is taken over all sections on $\mathcal{H}_1$. In general, the infimum is not necessarily attained  on
the whole space $W^{1,2}(\Sigma,\mathcal{L})$.

\begin{theorem}\label{T2}
Let $(\mathcal{L},\mathfrak{g})$ be a line bundle over a closed Riemann surface $(\Sigma,g)$, $\mathcal{D}$ be its bundle connection.
Suppose there exists a global unit frame $\zeta$ on $\Gamma(\mathcal{L})$. Let $\mathcal{J}_\rho$, $\mathcal{H}_0$ and $\mathcal{H}_1$ be defined as in (\ref{J-funct}), (\ref{H0}) and (\ref{H1}) respectively. Suppose that $\inf_{\sigma\in\mathcal{H}_1}\mathcal{J}_{8\pi}(\sigma)$ is not attained
on $\mathcal{H}_1$. Then we have the following two assertions:

$(i)$ If ${\rm dim}\mathcal{H}_0=0$, then there exist some $p_0\in\Sigma$ and  a Green section ${\rm G}_0=G_{p_0}\zeta$ with $$\Delta_{\mathcal{L}}{\rm G}_0=8\pi\le(\delta_{p_0}-\f{1}{|\Sigma|}\ri)\zeta$$
such that
\be
\inf_{\sigma\in\mathcal{H}_1}\mathcal{J}_{8\pi}(\sigma)=-8\pi-4\pi A_{p_0}-8\pi\log{\pi}-8\pi\log{h(p_0)}
+\frac{4\pi}{|\Sigma|}\int_{\Sigma}G_{p_0}dv_g,\label{lowerb-0}
\ee
where $A_{p_0}=\lim_{x\rightarrow p_0}(G_{p_0}(x)+4\log{d_g(x,p_0)})$ is a constant, $\Delta_{\mathcal{L}}=\mathcal{D}^\ast\mathcal{D}$
is the Laplace-Beltrami operator on line bundle,
$\mathcal{D}^\ast$ is the dual operator with respect to $\mathcal{D}$;

$(ii)$ If ${\rm dim}\mathcal{H}_0=1$ and $\{\tau_1\}$ is a normal basis of $\mathcal{H}_0$, then
there exist some $p\in\Sigma$ and  a Green section ${\rm G}=G_p\zeta $ with
\be
\le\{\begin{array}{lll}\Delta_{\mathcal{L}}{\rm G}=8\pi\le(\delta_p-\f{1}{|\Sigma|}\ri)\zeta-\lambda_{1}\tau_1
   \\[1.5ex]
   \lambda_{1}=8\pi\le(\langle\zeta,\tau_1\rangle(p)-\f{1}{|\Sigma|}\int_\Sigma\langle\zeta, \tau_1\rangle dv_g\ri)
\end{array}\ri.\label{green-1}
\ee
such that
\be
\inf_{\sigma\in\mathcal{H}_1}\mathcal{J}_{8\pi}(\sigma)=-8\pi-4\pi A_p-8\pi\log{\pi}-8\pi\log{h(p)}
+\frac{4\pi}{|\Sigma|}\int_{\Sigma}G_pdv_g,\label{lowerb}
\ee
where $A_{p}=\lim_{x\rightarrow p}(G_{p}(x)+4\log{d_g(x,p)})$ is a constant.
\end{theorem}

Indeed, Theorem \ref{T2} gives a sufficient condition such that $\mathcal{J}_{8\pi}$ has a minimizer on $\mathcal{H}_1$.
For the proof of Theorem \ref{T2}, we modify the method of blow-up
analysis used by Ding-Jost-Li-Wang \cite{DJLW} (see also \cite{LiLi,Wang,Yang-Zhu-Sci}).
 Several technical difficulties need to be overcome. Since $\mathcal{J}_\rho:W^{1,2}(\Sigma,\mathcal{L})\ra\mathbb{R}$
has no translation invariance, the method in \cite{DJLW} cannot be directly applied here. We overcome this difficulty by employing
a technique from \cite{Yang-Zhu-Sci}. As we shall see later, we need much more analysis on neck domains than in
\cite{DJLW,LiLi,Wang,Yang-Zhu-Sci}. Moreover, we need to construct a sequence of sections $\sigma_k$ satisfying
$\mathcal{J}_{8\pi}(\sigma_k)$ converges to $\inf_{\sigma\in\mathcal{H}_1}\mathcal{J}_{8\pi}(\sigma)$. This is also different from \cite{DJLW,LiLi,Wang}.

In the case $\zeta$ is a global unit frame on $\Gamma(\mathcal{L})$, if we further assume $\mathcal{D}\zeta=0$, then
we shall prove later that $\mathcal{H}_0={\rm span}\{\zeta\}$. Hence ${\rm dim}\mathcal{H}_0=1$, $\tau_1=|\Sigma|^{-1/2}\zeta$, and  (\ref{E-L}) would reduce to
\be\label{mean-eqn}\Delta_{\mathcal{L}}\sigma=\rho\le(\f{ he^{\langle\sigma,\zeta\rangle}}{\int_\Sigma he^{\langle\sigma,\zeta\rangle}dv_g}-
\f{1}{|\Sigma|}\ri)\zeta.\ee
Also we shall show  $\Delta_{\mathcal{L}}\sigma=(\Delta_gu)\zeta$ for any $\sigma=u\zeta$ in this case, and thus (\ref{mean-eqn})
is essentially the mean field equation (\ref{5}). As a consequence, our third result says

\begin{theorem}\label{Thm-2}
In addition to the assumptions in Theorem \ref{T1}, we further require that $\mathcal{D}\zeta=0$. Then we have the following
three assertions:\\
$(i)$ When $\rho<8\pi$, $(\ref{mean-eqn})$ has a solution;\\
$(ii)$ When $\rho=8\pi$, if the Ding-Jost-Li-Wang condition is satisfied, then  $(\ref{mean-eqn})$ has a solution;\\
$(iii)$ When $\rho\in(8k\pi,8(k+1)\pi)$, $\forall k\in\mathbb{N}^\ast$, $(\ref{mean-eqn})$ has a solution.
\end{theorem}

Here the {\it Ding-Jost-Li-Wang condition} means certain assumptions on $h$ and $(\Sigma,g)$ proposed by Ding-Jost-Li-Wang
(\cite{DJLW}, Theorem 1.2). For the proof of Theorem \ref{Thm-2}, we only need to transform the problem from line bundle setting to
Riemann surface setting, and then use the existing results of Ding-Jost-Li-Wang \cite{DJLW} and Djadli \cite{Djadli}.

\bigskip
The remaining part of this paper is organized as follows:  In section 2, we give some preliminaries;
In Section 3, by a direct method of variation, we prove Theorem \ref{T1};
In Section 4, we use the method of blow-up analysis to prove Theorem \ref{T2};
In Section 5, we transform Theorem \ref{Thm-2} to results of Ding-Jost-Li-Wang \cite{DJLW} and Djadli \cite{Djadli}.
Throughout this paper, we do not distinguish sequence and subsequence, often denote various constants by the same $C$, and represent
metrics on different bundles by the same $|\cdot|$. We always denote a geodesic ball centered at $x\in\Sigma$ with radius $r$ by
$B_r(x)=\{q\in\Sigma:d_g(q,x)<r\}$, and a Euclidean ball centered at $0$ with radius $r$ by $\mathbb{B}_r=\{z\in\mathbb{R}^2: |z|<r\}$.

\section{Preliminaries}
In this section, we present some useful facts which will be used frequently in the
following.
First, a distributional decomposition of $\Delta_{\mathcal{L}}$ will be given.
Next, we show the Poincar\'e inequality holds on the space $\mathcal{H}_1$.
Finally, we derive several properties of the spaces
$\mathcal{H}_0$ and $\mathcal{H}_1$, which are defined by (\ref{H0}) and (\ref{H1}) respectively.

To write $\Delta_{\mathcal{L}}$ explicitly, we
recall
$\zeta$ is a global unit frame on $W^{1,2}(\Sigma,\mathcal{L})$ and $\mathcal{D}\zeta=\omega\otimes\zeta$. Let $d^\ast: W^{1,2}(\Sigma,T^\ast\Sigma)\ra W^{1,2}(\Sigma)$ be defined by
$$\int_\Sigma (d^\ast\varsigma) u dv_g=\int_\Sigma \langle\varsigma,du\rangle dv_g$$
for all $\varsigma\in W^{1,2}(\Sigma,T^\ast\Sigma)$ and all $u\in W^{1,2}(\Sigma)$.  In a coordinate system $\{x^\alpha\}_{\alpha=1}^2$, if $\varsigma=\varsigma_\alpha dx^\alpha$, then one can easily compute
\be\label{dstar}d^\ast\varsigma=-\f{1}{\sqrt{g}}\f{\p}{\p x^\alpha}\le(g^{\alpha\beta}\sqrt{g}\varsigma_\beta\ri),\ee
where $(g^{\alpha\beta})$ is the inverse matrix of $(g_{\alpha\beta})$ and $\sqrt{g}=\sqrt{{\rm det}(g_{\alpha\beta})}$. Obviously the
Laplace-Beltrami operator $\Delta_g$ is locally written as
\be\label{Laplace}\Delta_gu=d^\ast du=-\f{1}{\sqrt{g}}\f{\p}{\p x^\alpha}\le(g^{\alpha\beta}\sqrt{g}\f{\p u}{\p x^\beta}\ri).\ee

A straightforward calculation gives

\begin{lemma}\label{delta-L}
For any $\sigma=u\zeta\in W^{1,2}(\Sigma,\mathcal{L})$, there is a distributional decomposition
$$\Delta_{\mathcal{L}}\sigma=\le(\Delta_gu+(|\omega|^2+d^\ast\omega)u\ri)\zeta.$$
\end{lemma}
\proof Let $\sigma=u\zeta\in W^{1,2}(\Sigma,\mathcal{L})$. Take any $\phi=v\zeta\in \Gamma(\mathcal{L})$. Note that $|\zeta|=1$.
In view of (\ref{dstar}) and (\ref{Laplace}), one has
\bna
\int_\Sigma \langle\sigma,\Delta_\mathcal{L}\phi\rangle dv_g&=&
\int_\Sigma\langle\mathcal{D}\sigma,\mathcal{D}\phi\rangle dv_g\\
&=&\int_\Sigma\langle du\otimes\zeta+u\omega\otimes\zeta,dv\otimes\zeta+v\omega\otimes\zeta\rangle dv_g\\
&=&\int_\Sigma \le(d^\ast du+(|\omega|^2+d^\ast\omega)u\ri)vdv_g\\
&=&\int_\Sigma\langle (\Delta_gu+(|\omega|^2+d^\ast\omega)u)\zeta,\phi\rangle dv_g.
\ena
This is exactly the desired result.$\hfill\Box$
\begin{lemma}\label{finite-dimension}
   Any section in $\mathcal{H}_0$ is smooth on $\Sigma$. Moreover, 
    $\mathcal{H}_0$ is a finite dimensional linear space, in particular ${\rm dim}\mathcal{H}_0\leq 1$.
    \end{lemma}
    \proof
    Taking any section $\tau$ in $\mathcal{H}_0$, since for all $\varphi\in\Gamma(\mathcal{L})$,
 $$0=\int_\Sigma\langle\mathcal{D}\tau,\mathcal{D}\varphi\rangle dv_g=\int_\Sigma\langle\tau,\mathcal{D}^\ast\mathcal{D}
 \varphi\rangle dv_g=\int_\Sigma\langle\tau,\Delta_{\mathcal{L}}\varphi\rangle dv_g,$$
 we know that $\Delta_{\mathcal{L}}\tau=0$ in the distributional sense. Let $\tau=\phi\zeta$. By Lemma \ref{delta-L},
 $\phi$ is a distributional solution of
 \be\label{phi}\Delta_g\phi+(|\omega|^2+d^\ast\omega)\phi=0.\ee
 Noting $\phi\in W^{1,2}(\Sigma)$, we have by repeatedly applying elliptic estimate (\cite{Aubin1982}, Theorem 3.54)
 to (\ref{phi}) that $\phi\in W^{\ell,2}(\Sigma)$ for all $\ell\geq 2$.
 Then the Sobolev embedding theorem leads to $\phi\in C^\infty(\Sigma)$. Hence,
 we obtain $\tau\in\Gamma(\mathcal{L})$, which confirms the first assertion.
 
   For any $\sigma\in\mathcal{H}_0$, there is a unique $u\in C^\infty(\Sigma)$ such that $\sigma=u\zeta$ and
    $\mathcal{D}\sigma=0$. This is equivalent to the equation
    \be\label{zero}du+u\omega=0.\ee
    Obviously $u\equiv 0$ is a solution of (\ref{zero}). If (\ref{zero}) has a unique solution $u\equiv 0$, 
     then ${\rm dim}\mathcal{H}_0= 0$. In the following, we consider the case $\mathcal{H}_0\not=\{0\}$. We first {\it claim} that if $u$ is a solution of (\ref{zero}) and
    $u(x_0)=0$ for some $x_0\in \Sigma$, then
    $u\equiv 0$ on $\Sigma$. To see this, for any $x\in\Sigma$ and $x\not=x_0$, we take a geodesic line $\gamma:[0,1]\ra \Sigma$ satisfying
    $\gamma(0)=x_0$ and $\gamma(1)=x$. Set $y(t)=u(\gamma(t))$ and $a(t)=\omega(\gamma(t))(\gamma^\prime(t))$ for $t\in [0,1]$, where $\gamma^\prime(t)$ denotes the tangent vector field along the geodesic line $\gamma(t)$. It follows from (\ref{zero}) that
    \be\label{ode}\le\{\begin{array}{lll}
    \f{d}{dt}y(t)=-a(t)y(t)\\[1.5ex]
    y(0)=0.
    \end{array}\ri.\ee
    Applying the existence and uniqueness theorem of the linear ODE to (\ref{ode}), we obtain $y(t)\equiv 0$ for all $t\in[0,1]$. Thus 
    $u(x)=u(\gamma(1))=y(1)=0$. Since $x$ is an arbitrary point on $\Sigma$, $u$ must be identically zero, and our first claim follows.
    In other words, for any $\sigma\in\mathcal{H}_0$, either $\sigma(x)\equiv 0$ for all $x\in\Sigma$, or $\sigma(x)\not= 0$ for all $x\in\Sigma$.
    Our second {\it claim} is that for any two nonzero sections $\sigma=u\zeta\in\mathcal{H}_0$ and $\tau=v\zeta\in\mathcal{H}_0$, there exists
    a constant $c$ such that $\sigma=c\tau$. Obviously the first claim implies that $u(x)\not=0$ and $v(x)\not=0$ for all $x\in\Sigma$. This together with (\ref{zero}) gives
    \be\label{equal}\f{du}{u}=\f{dv}{v}\quad{\rm on}\quad \Sigma.\ee
    Clearly (\ref{equal}) leads to  $\log|u/v|\equiv C$ for some constant $C$, and there exists some constant $c$ such that $u(x)=cv(x)$ for all $x\in\Sigma$.
    This confirms the second claim, and implies ${\rm dim}\mathcal{H}_0=1$.
    $\hfill\Box$
    \\
 
By Lemma \ref{finite-dimension}, ${\rm dim}\mathcal{H}_0\leq 1$. If ${\rm dim}\mathcal{H}_0=0$, then
$\mathcal{H}_0=\{0\}$; If ${\rm dim}\mathcal{H}_0=1$, then we may assume
$\mathcal{H}_0={\rm span}\{\tau_1\}$, where $\{\tau_1\}$ is a normal basis on $\mathcal{H}_0$.
We next show that the Poincar\'e inequality holds on the space $\mathcal{H}_1$.
\begin{lemma}\label{Poincare}
There exists a constant $C$ depending only on $(\Sigma,g)$ and $(\mathcal{L},\mathfrak{g})$ such that
$$\int_\Sigma|\sigma|^2dv_g\leq C\int_\Sigma|\mathcal{D}\sigma|^2dv_g$$
for all $\sigma\in\mathcal{H}_1$.
\end{lemma}
\proof
Suppose not. For any $k\in\mathbb{N}^\ast$, there would be a section $\sigma_k\in\mathcal{H}_1$ satisfying
\be\label{sk1}\int_\Sigma|\sigma_k|^2dv_g=1\ee
and
\be\label{der-0}\int_\Sigma|\mathcal{D}\sigma_k|^2dv_g\leq \f{1}{k}.\ee
Thus $\{\sigma_k\}$ is bounded in $W^{1,2}(\Sigma,\mathcal{L})$. With no loss of generality, we find some $\sigma_0\in
W^{1,2}(\Sigma,\mathcal{L})$ such that $\sigma_k$ converges to $\sigma_0$ weakly in $W^{1,2}(\Sigma,\mathcal{L})$,
strongly in $L^2(\Sigma,\mathcal{L})$, and almost everywhere in $\Sigma$. This implies  $\sigma_0\in\mathcal{H}_1$. Moreover (\ref{sk1}) and (\ref{der-0})
lead to
\be\label{s0}\int_\Sigma|\sigma_0|^2dv_g=1\ee
and
\be\label{ds0}\int_\Sigma|\mathcal{D}\sigma_0|^2dv_g=0.\ee
By (\ref{ds0}), we have $\sigma_0\in\mathcal{H}_0$. Hence $\sigma_0\in\mathcal{H}_0\cap\mathcal{H}_1=\{0\}$, which  contradicts (\ref{s0}).
In this way, we get the desired result. $\hfill\Box$

\section{The subcritical case}
In this section, we prove Theorem \ref{T1}. Since the proof in the case ${\rm dim}\mathcal{H}_0=0$ is a minor modification of that
in the case ${\rm dim}\mathcal{H}_0=1$, we only prove the theorem under the assumption ${\rm dim}\mathcal{H}_0=1$, but leave
the case ${\rm dim}\mathcal{H}_0=0$ to interested readers.
To begin with,  we pay attention to the Trudinger-Moser inequality on line bundle. Though a more general inequality was already established in \cite{LiLiuYang},
 we will give a simple proof in our setting. Also this can be compared with \cite{Yang-06}.

 \begin{lemma}\label{small}
 There exists some $\alpha_0>0$ such that
 \be\label{alp0}\sup_{\sigma\in\mathcal{H}_1,\int_\Sigma|\mathcal{D}\sigma|^2dv_g\leq 1}
 \int_\Sigma e^{\alpha_0|\sigma|^2}dv_g<+\infty.\ee
 \end{lemma}
 \proof
 Since $\zeta$ is a global unit frame on $W^{1,2}(\Sigma,\mathcal{L})$, for any $\sigma\in W^{1,2}(\Sigma,\mathcal{L})$, there must be
 a unique function $u\in W^{1,2}(\Sigma)$ such that $\sigma=u\zeta$ and $|\sigma|=|u|$.
 If $\sigma\in\mathcal{H}_1$ satisfies $\int_\Sigma|\mathcal{D}\sigma|^2dv_g\leq 1$, then Lemma \ref{Poincare} implies
 $\int_\Sigma|\sigma|^2dv_g\leq C$, and thus  $\int_\Sigma u^2dv_g\leq C$. This together with the formula (\ref{Lebn}) leads to
 \bna
 \int_\Sigma|du|^2dv_g&\leq&2\int_\Sigma|\mathcal{D}\sigma|^2dv_g+2\int_\Sigma u^2|\omega|^2dv_g\\
 &\leq&C.
 \ena
 It then follows that $\|u\|_{W^{1,2}(\Sigma)}\leq C_0$ for some constant $C_0$. By the classical Trudinger-Moser inequality \cite{Moser71,Fontana}, we have
 $$\int_\Sigma e^{\f{4\pi}{C_0^2}u^2}dv_g\leq S_{4\pi}=\sup_{\|u\|_{W^{1,2}(\Sigma)}\leq 1}\int_\Sigma e^{4\pi u^2}dv_g\leq
 C\sup_{\int_\Sigma|d u|^2dv_g\leq 1,\int_\Sigma udv_g=0}\int_\Sigma e^{4\pi u^2}dv_g.$$
 Equivalently
 $$\int_\Sigma e^{\f{4\pi}{C_0^2}|\sigma|^2}dv_g\leq C.$$
 One can see that any $\alpha_0\in(0,4\pi/C_0^2]$ satisfies (\ref{alp0}), as desired. $\hfill\Box$\\

  Define $\alpha_\ast=\sup\{\alpha:\sup_{\sigma\in\mathcal{H}_1,\int_\Sigma|\mathcal{D}\sigma|^2dv_g\leq 1}
  \int_\Sigma e^{\alpha|\sigma|^2}dv_g<+\infty\}$.
  The following embedding is analogous to the usual one for $W^{1,2}(\Sigma)$.

 \begin{lemma}\label{arbitr}
 For any $\sigma\in W^{1,2}(\Sigma,\mathcal{L})$, there holds
 $$\int_\Sigma e^{p|\sigma|^2}dv_g<+\infty,\quad\forall p>1.$$
 \end{lemma}
 \proof
 Let $\sigma\in W^{1,2}(\Sigma,\mathcal{L})$ be fixed. There exists a unique function $u$ such that $\sigma=u\zeta$. In view of
 (\ref{Lebn}), we have $u\in W^{1,2}(\Sigma)$. Then the results of \cite{Moser71,Fontana} imply
 $$\int_\Sigma e^{p|\sigma|^2}dv_g=\int_\Sigma e^{pu^2}dv_g<+\infty$$
 for any $p>1$.
 $\hfill\Box$\\

 Now we are in a position to estimate $\alpha_\ast$.

 \begin{lemma}\label{best}
 $\alpha_\ast=4\pi$.
 \end{lemma}
 \proof {\it Step} 1. There holds $\alpha_\ast\geq 4\pi$.

 Suppose not. In view of Lemma \ref{small},
 $0<\alpha_\ast<4\pi$. According to the definition of $\alpha_\ast$, for any $j\in\mathbb{N}^\ast$, there would exist some $\sigma_j\in \mathcal{H}_1$ such that
 $\int_\Sigma|\mathcal{D}\sigma_j|^2dv_g\leq 1$ and
 \be\label{infty}\int_\Sigma e^{(\alpha_\ast+j^{-1})|\sigma_j|^2}dv_g\ra+\infty\quad{\rm as}\quad j\ra\infty.\ee
 By Lemma \ref{Poincare}, $\{\sigma_j\}$ is bounded in $W^{1,2}(\Sigma,\mathcal{L})$. With no loss of generality,
 we assume $\sigma_j$ converges to $\sigma^\ast$ weakly in $W^{1,2}(\Sigma,\mathcal{L})$, strongly in $L^2(\Sigma,\mathcal{L})$,
 and almost everywhere in $\Sigma$. Clearly $\sigma^\ast\in\mathcal{H}_1$. Moreover, we {\it claim} that $\sigma^\ast=0$.
 For otherwise, we have $\int_\Sigma|\mathcal{D}\sigma^\ast|^2dv_g>0$ and
 \bna
 \int_\Sigma|\mathcal{D}(\sigma_j-\sigma^\ast)|^2dv_g&=&\int_\Sigma|\mathcal{D}\sigma_j|^2dv_g-
 \int_\Sigma|\mathcal{D}\sigma^\ast|^2dv_g+o_j(1)\nonumber
 \\&\leq&1-\f{1}{2}\int_\Sigma|\mathcal{D}\sigma^\ast|^2dv_g,
 \ena
 provided that $j$ is chosen sufficiently large. Then for any $\epsilon>0$, the Young inequality leads to
 \bna
 |\sigma_j|^2&\leq&(1+\epsilon)|\sigma_j-\sigma^\ast|^2+(1+\epsilon^{-1})|\sigma^\ast|^2\\
 &\leq&(1+\epsilon)\le(1-\f{1}{2}\int_\Sigma|\mathcal{D}\sigma^\ast|^2dv_g\ri)\f{|\sigma_j-\sigma^\ast|^2}
 {\int_\Sigma|\mathcal{D}(\sigma_j-\sigma^\ast)|^2dv_g}+(1+\epsilon^{-1})|\sigma^\ast|^2,
 \ena
 if $j$ is large enough. Fixing $\epsilon$ with $0<\epsilon<\f{1}{8}\int_\Sigma|\mathcal{D}\sigma^\ast|^2dv_g$, we obtain
 \bna(\alpha_\ast+j^{-1})|\sigma_j|^2&\leq& (\alpha_\ast+j^{-1})(1+\epsilon)\le(1-\f{1}{2}\int_\Sigma|\mathcal{D}\sigma^\ast|^2dv_g\ri)
 \f{|\sigma_j-\sigma^\ast|^2}
 {\int_\Sigma|\mathcal{D}(\sigma_j-\sigma^\ast)|^2dv_g}+C(\epsilon)|\sigma^\ast|^2\\
 &\leq& \le(\alpha_\ast-\f{1}{4}\int_\Sigma|\mathcal{D}\sigma^\ast|^2dv_g\ri)\f{|\sigma_j-\sigma^\ast|^2}
 {\int_\Sigma|\mathcal{D}(\sigma_j-\sigma^\ast)|^2dv_g}+C(\epsilon)|\sigma^\ast|^2\ena
 for sufficiently large $j$. By the definition of $\alpha^\ast$, Lemma \ref{arbitr} and the H\"older inequality, we have
 $$\int_\Sigma e^{(\alpha_\ast+j^{-1})|\sigma_j|^2}dv_g\leq C.$$
 This contradicts (\ref{infty}), and confirms our claim $\sigma^\ast=0$.

 Write $\sigma_j=u_j\zeta$. From the above, we know that $u_j$ converges to $0$ weakly in $W^{1,2}(\Sigma)$, strongly in $L^2(\Sigma)$,
  and almost everywhere in $\Sigma$. On one hand, there holds for any $\delta>0$,
 \bna
 \int_\Sigma|\mathcal{D}\sigma_j|^2dv_g&=&\int_\Sigma|du_j\otimes\zeta+u_j\omega\otimes\zeta|^2dv_g\\
 &\geq&(1-\delta)\int_\Sigma|du_j|^2dv_g+o_j(1).
 \ena
 It follows that for sufficiently large $j$,
 \bna\label{grad}\int_\Sigma|du_j|^2dv_g\leq \f{1}{1-2\delta}.\ena
 Choosing $\delta>0$ small enough such that $(\alpha_\ast+j^{-1})/(1-2\delta)<4\pi-\delta$, and noting that $u_j$ converges to $0$ weakly in
 $W^{1,2}(\Sigma)$, strongly in $L^2(\Sigma)$,
 we obtain by the classical Trudinger-Moser inequality \cite{Moser71,Fontana},
 $$\int_\Sigma e^{(\alpha_\ast+j^{-1})|\sigma_j|^2}dv_g=\int_\Sigma e^{(\alpha_\ast+j^{-1})|u_j|^2}dv_g\leq C.$$
 This contradicts (\ref{infty}), and ends the first step.

 {\it Step} 2. There holds $\alpha_\ast\leq 4\pi$.

  Take $z\in \Sigma$ and $\delta>0$. For $k\in\mathbb{N}^\ast$, we define
    $$u_k=\left\{\begin{array}{lll}
                     -\sqrt{\f{1}{4\pi}\log k}
                     &{\rm in}&
                     B_{\delta/\sqrt{k}}(z)\\[1.5ex]
                      \f{1}{\sqrt{\pi\log k}}\log\f{r}{\delta}
                      &{\rm in}& B_\delta(z)\setminus B_{\delta/\sqrt{k}}(z)
                      \\[1.5ex]
                      0 &{\rm in}& \Sigma\setminus
                      B_{\delta}(z)
                   \end{array}\right.
    $$
    and
    $$\sigma_k=u_k\zeta-\le(\int_\Sigma\langle u_k\zeta,\tau_1\rangle dv_g\ri)\tau_1,$$
    where $r=r(x)=d_g(x,z)$ and  $\{\tau_1\}$ is a normal basis on $\mathcal{H}_0$.
    Obviously $\sigma_k\in \mathcal{H}_1$. In view of Lemma \ref{finite-dimension}, there exists some constant
    $C=\max_{x\in\Sigma}|\tau_1(x) | .$
    For any $\alpha>4\pi$, one can check that if $\delta>0$ is chosen sufficiently
    small, then
    \bna\sup_{\sigma\in\mathcal{H}_1,\int_\Sigma|\mathcal{D}\sigma|^2dv_g\leq 1}\int_\Sigma e^{\alpha|\sigma|^2}dv_g&\geq&
    \int_\Sigma \exp\le\{\alpha\f{|\sigma_k|^2}{\int_\Sigma|\mathcal{D}\sigma_k|^2dv_g}\ri\}dv_g\\
    &\geq&\int_{B_{\delta/\sqrt{k}}(z)}\exp\le\{\f{\alpha}{4\pi}(1+o_\delta(1)+o_k(1))\log k\ri\}dv_g\\
    &=&\pi\delta^2(1+o_\delta(1))k^{\f{\alpha}{4\pi}(1+o_\delta(1)+o_k(1))-1}\\
    &\ra& +\infty\quad{\rm as}\quad k\ra\infty.\ena
    This ends the second step.\\

    The lemma follows immediately from the above two steps. $\hfill\Box$\\

   {\it Proof of Theorem \ref{T1}}. The proof is based on a variational  method. Fixing $\rho<8\pi$, we set
   $$\Lambda_\rho=\inf_{\sigma\in\mathcal{H}_1}\mathcal{J}_\rho(\sigma).$$
   According to Lemma \ref{best}, $\Lambda_\rho>-\infty$. Take a minimizing sequence of sections $\sigma_k\in\mathcal{H}_1$
   such that $\mathcal{J}_\rho(\sigma_k)\ra \Lambda_\rho$ as $k\ra\infty$. By Lemmas \ref{Poincare} and \ref{best}, we have for any
   $\alpha<4\pi$,
   \bna
   \Lambda_\rho+o_k(1)&=&\f{1}{2}\int_\Sigma|\mathcal{D}\sigma_k|^2dv_g+\f{\rho}{|\Sigma|}\int_\Sigma \langle\sigma_k,\zeta\rangle dv_g
   -\rho\log\int_\Sigma h e^{\langle\sigma_k,\zeta\rangle} dv_g\\
   &\geq&\le(\f{1}{2}-\f{\rho}{4\alpha}-\epsilon\ri)\int_\Sigma|\mathcal{D}\sigma_k|^2dv_g-C(\alpha,\epsilon).
   \ena
   Choosing $\alpha$ satisfying $\rho/2<\alpha<4\pi$ and sufficiently small $\epsilon>0$, we conclude that $\{\sigma_k\}$ is bounded in $W^{1,2}(\Sigma,\mathcal{L})$.
   Hence, there exists some $\sigma_0$ such that $\sigma_k$ converges to $\sigma_0$ weakly in $W^{1,2}(\Sigma,\mathcal{L})$, strongly in
   $L^2(\Sigma,\mathcal{L})$, and almost everywhere in $\Sigma$. It follows that $\sigma_0\in\mathcal{H}_1$,
   $$\int_\Sigma|\mathcal{D}\sigma_0|^2dv_g\leq \limsup_{k\ra\infty}\int_\Sigma|\mathcal{D}\sigma_k|^2dv_g,$$
   that
   $$\int_\Sigma \langle\sigma_0,\zeta\rangle dv_g=\lim_{k\ra\infty}\int_\Sigma \langle\sigma_k,\zeta\rangle dv_g,$$
   and that
   $$\int_\Sigma h e^{\langle\sigma_0,\zeta\rangle} dv_g=\lim_{k\ra\infty}\int_\Sigma h e^{\langle\sigma_k,\zeta\rangle} dv_g.$$
   Summarizing the above three estimates, we conclude
   $$\Lambda_\rho\leq\mathcal{J}_\rho(\sigma_0)\leq \lim_{k\ra\infty}\mathcal{J}_\rho(\sigma_k)=\Lambda_\rho.$$
   Therefore $\sigma_0$ is a minimizer of $\mathcal{J}_\rho$ on $\mathcal{H}_1$. This completes the proof of the theorem. $\hfill\Box$

   \section{Critical case}
   In this section, we shall prove Theorem \ref{T2} by using the blow-up analysis. Since the functional $\mathcal{J}_\rho$
   has no translation invariance, and the maximum principle is not available in the line bundle setting, our argument is quite different from that of \cite{DJLW,LiLi,Wang,Yang-Zhu-Sci}. Since the consideration of the case ${\rm dim}\mathcal{H}_0=0$ is almost the same
   as that of the case ${\rm dim}\mathcal{H}_0=1$, we only consider the second case, and leave
   the first case to interested readers.

   \subsection{Maximizers for subcritical functionals}

   As Ding-Jost-Li-wang did in \cite{DJLW}, we shall analyze maximizers for subcritical functionals.
   For any positive integer $k$, we denote $\rho_k=8\pi-1/k$. By Theorem \ref{T1}, there exists a section sequence
   $\sigma_k\in\mathcal{H}_1$ satisfying $\mathcal{J}_{\rho_k}(\sigma_k)=\inf_{\sigma\in\mathcal{H}_1}\mathcal{J}_{\rho_k}(\sigma)$.
   Clearly $\sigma_k$ satisfies the Euler-Lagrange equation
   \be\label{E-Lk}\le\{\begin{array}{lll}\Delta_{\mathcal{L}}\sigma_k=\rho_k\le(\f{ he^{\langle\sigma_k,\zeta\rangle}}{\mu_k}-
\f{1}{|\Sigma|}\ri)\zeta-\lambda_{1,k}\tau_1\\[1.5ex]
\mu_k=\int_\Sigma he^{\langle\sigma_k,\zeta\rangle}dv_g
\\[1.5ex]
\lambda_{1,k}=\rho_k\int_\Sigma\le(\f{ he^{\langle\sigma_k,\zeta\rangle}}{\mu_k}
-\f{1}{|\Sigma|}\ri) \langle\zeta,\tau_1\rangle dv_g,
\end{array}\ri.\ee
where $\Delta_\mathcal{L}=\mathcal{D}^\ast\mathcal{D}$ and $\{\tau_1\}$ is a normal basis on $\mathcal{H}_0$.

In view of the importance of the coefficients $\lambda_{1,k}$ and $\mu_k$,  we have the following
\begin{lemma}\label{muk}
There exists a constant $C$ depending only on $(\Sigma,g)$ and $(\mathcal{L},\mathfrak{g})$ such that $|\lambda_{1,k}|\leq C$ 
for all $k\in\mathbb{N}^\ast$. Moreover, there holds
$$\liminf_{k\ra\infty}\mu_k>0.$$
\end{lemma}
\proof By Lemma \ref{finite-dimension}, we assume $\mathcal{H}_0={\rm span}\{\tau_1\}$ and $\tau_1\in \Gamma(\mathcal{L})$.
This together with the facts $|\zeta(x)|=1$, $h(x)>0$ for all $x\in\Sigma$, and $0<\rho_k\leq 8\pi$ implies
$|\lambda_{1,k}|\leq C$ for all $k\in\mathbb{N}^\ast$.

For the second assertion, we suppose not. Then up to a subsequence, $\mu_k\ra 0$ as $k\ra\infty$. One easily sees
from the choice of $\sigma_k$ that
\bea\nonumber
\mathcal{J}_{\rho_k}(\sigma_k)&=&\f{1}{2}\int_\Sigma|\mathcal{D}\sigma_k|^2dv_g+\f{\rho_k}{|\Sigma|}\int_\Sigma
\langle\sigma_k,\zeta\rangle dv_g-\rho_k\log\mu_k\\\nonumber
&=&\inf_{\sigma\in\mathcal{H}_1}\mathcal{J}_{\rho_k}(\sigma)\\\nonumber
&\leq&\mathcal{J}_{\rho_k}(0)\\\label{les-1}
&\leq&8\pi\le|\log\int_\Sigma h dv_g\ri|,
\eea
and from the Poincar\'e inequality that
\bea\nonumber
\le|\int_\Sigma\langle\sigma_k,\zeta\rangle dv_g\ri|&\leq&\int_\Sigma|\sigma_k|dv_g\\\nonumber
&\leq&C\le(\int_\Sigma|\mathcal{D}\sigma_k|^2dv_g\ri)^{1/2}\\\label{les-2}
&\leq&\f{1}{4}\int_\Sigma|\mathcal{D}\sigma_k|^2dv_g+C.
\eea
Combining (\ref{les-1}) and (\ref{les-2}), one has
$\rho_k\log\mu_k\geq -C$ for some positive constant $C$. This contradicts $\mu_k$ converges to $0$, and ends the proof of the lemma.
$\hfill\Box$

\begin{lemma}\label{lim-8pi} There holds
\be\label{limit}\lim_{k\ra\infty}\inf_{\sigma\in\mathcal{H}_1}\mathcal{J}_{\rho_k}(\sigma)=
\inf_{\sigma\in\mathcal{H}_1}\mathcal{J}_{8\pi}(\sigma).\ee
\end{lemma}
\proof Since the proof is completely analogous to that of (\cite{Yang-Zhu-Sci}, Lemma 2.3),
we omit the details here.$\hfill\Box$

\begin{lemma}\label{no-blow1}
If $\mu_k$ is bounded, then $\sigma_k$ is bounded in $W^{1,2}(\Sigma,\mathcal{L})$ and $\mathcal{J}_{8\pi}$ has a minimizer
on $\mathcal{H}_1$.
\end{lemma}
\proof If $\mu_k$ is bounded, then (\ref{les-1}) gives
$$\int_\Sigma|\mathcal{D}\sigma_k|^2dv_g\leq C.$$
This implies that $\sigma_k$ is bounded in $W^{1,2}(\Sigma,\mathcal{L})$, and thus
$\sigma_k$ converges to some $\sigma_0$ weakly in $W^{1,2}(\Sigma,\mathcal{L})$, strongly in $L^2(\Sigma,\mathcal{L})$,
and almost everywhere in $\Sigma$. This together with (\ref{limit}) leads to that $\sigma_0\in\mathcal{H}_1$ and
$\mathcal{J}_{8\pi}(\sigma_0)=\inf_{\sigma\in\mathcal{H}_1}\mathcal{J}_{8\pi}(\sigma)$. $\hfill\Box$\\

From now on in this section, we write $\sigma_k=u_k\zeta$ and
\be\label{ck}c_k=\max_\Sigma u_k=u_k(x_k).\ee

\begin{lemma}\label{c-bdd}
If $c_k$ is bounded from above, then $\mathcal{J}_{8\pi}$ has a minimizer on $\mathcal{H}_1$.
\end{lemma}
\proof In view of Lemma \ref{muk}, multiplying both sides of (\ref{E-Lk}) by $\sigma_k$ and integrating by parts, one estimates by
the H\"older inequality
\bna
\int_\Sigma|\mathcal{D}\sigma_k|^2dv_g&=&\rho_k\int_\Sigma\le(\f{he^{u_k}}{\mu_k}-\f{1}{|\Sigma|}\ri)\langle\zeta,\sigma_k
\rangle dv_g-\lambda_{1,k}\int_\Sigma\langle\tau_1,\sigma_k\rangle dv_g\\
&\leq&\f{\rho_k}{\mu_k}\le(\int_\Sigma h^2e^{2u_k}dv_g\ri)^{1/2}\le(\int_\Sigma|\sigma_k|^2dv_g\ri)^{1/2}+C\le(\int_\Sigma|\sigma_k|^2dv_g\ri)^{1/2}\\
&\leq& C\le(\int_\Sigma|\mathcal{D}\sigma_k|^2dv_g\ri)^{1/2}.
\ena
This implies that $\sigma_k$ is bounded in $W^{1,2}(\Sigma,\mathcal{L})$. Arguing similarly as in the proof of Lemma \ref{no-blow1},
one concludes the lemma. $\hfill\Box$\\

If $\mathcal{J}_{8\pi}$ has no minimizer on $\mathcal{H}_1$, then it follows from Lemmas \ref{no-blow1} and \ref{c-bdd} that
\be\label{blowup}\mu_k\ra+\infty,\quad c_k\ra+\infty.\ee
Let $x_k$ be as in (\ref{ck}). With no loss of generality, we assume
\bna\label{concen-p}x_k\ra p\quad {\rm as}\quad k\ra\infty.\ena
 Take a sequence of isothermal coordinate systems
$\left(U, \phi_{k} ;\left\{y_{1}, y_{2}\right\}\right)$ near $x_k$ with $\phi_{k}\left(x_{k}\right)=0$
such that $\phi_k^{-1}(\mathbb{B}_{2\delta})\subset U$ and the metric $g$ is represented by
\bna\label{isotherm}g_{k}(y)=\exp ({f_k(y)})(d y^2_1+d y^2_2),\ena where
$f_{k}: \phi_{k}(U) \rightarrow \mathbb{R}$ is a smooth function satisfying
$f_{k}(0)=0$, $|d^jf_k|\leq C$, $j=1,2$,
\be\label{dist}\exp(f_k(y))=1+O(|y|)=1+O(d_g(\phi_k^{-1}(y),x_k)),\ee
$C^{-1}|y|\leq d_g(\phi_k^{-1}(y),x_k)\leq C|y|$ for some constant $C$ independent of $k$. Moreover
\bna\label{iso-trans}|du|_{g_k}^2=\exp(-f_k)|du|_{\mathbb{R}^2}^2,\quad \Delta_{g_k}=-\exp(-f_k)\Delta_{\mathbb{R}^2}.\ena
Indeed, such a sequence of isothermal coordinates exist. For its explicit proof, we refer the readers to (\cite{Yang2020}, Lemma 3.1).
\subsection{Blow-up analysis}
To analyze the asymptotic behavior of $\sigma_k$, we set
\be\label{scale}r_k=\f{\sqrt{\mu_k}}{\sqrt{\rho_k h(p)}}e^{-\f{c_k}{2}}.\ee
Recalling (\ref{blowup}), we have an analog of (\cite{Yang-Zhu-Sci}, Lemma 2.7), namely
\begin{lemma}\label{r-0}
For any $\gamma<1/2$, there holds $r_k^2e^{\gamma c_k}\ra 0$ as $k\ra\infty$. In particular for any $q>0$, there holds $r_kc_k^q\ra 0$ as $k\ra\infty$.
\end{lemma}
\proof Recall $h>0$ on $\Sigma$. Multiplying both sides of (\ref{E-Lk}) by $\sigma_k$, integrating by parts and
using the Poincar\'e inequality (see Lemma \ref{Poincare} above),  we have
\bea\nonumber
\int_\Sigma|\mathcal{D}\sigma_k|^2dv_g&=&\rho_k\int_\Sigma\le(\f{he^{u_k}}{\mu_k}-\f{1}{|\Sigma|}\ri)\langle\sigma_k,\zeta\rangle dv_g
-\lambda_{1,k}\int_\Sigma\langle\tau_1,\sigma_k\rangle dv_g\\\nonumber
&\leq&8\pi c_k\int_\Sigma\f{he^{u_k}}{\mu_k}dv_g+C\int_\Sigma|\sigma_k|dv_g\\\label{gr-0}
&\leq&8\pi c_k+C\le(\int_\Sigma|\mathcal{D}\sigma_k|^2dv_g\ri)^{1/2}.
\eea
Let $0<\epsilon<1$ be a fixed number to be determined later. We conclude from (\ref{gr-0}) that for all sufficiently large $k$,
\be\label{gr-2}\int_\Sigma|\mathcal{D}\sigma_k|^2dv_g\leq 8\pi(1+\epsilon)c_k.\ee
This together with Lemma \ref{best} gives
\bna
\int_\Sigma he^{\langle\sigma_k,\zeta\rangle}dv_g&\leq& (\max_\Sigma h)\int_\Sigma \exp\le\{
\f{4\pi(1-\epsilon)|\sigma_k|^2}{\int_\Sigma|\mathcal{D}\sigma_k|^2dv_g}+\f{\int_\Sigma|\mathcal{D}\sigma_k|^2dv_g}
{16(1-\epsilon)\pi}\ri\}dv_g\\
&\leq& C(\epsilon) \exp\le\{\f{1+\epsilon}{2(1-\epsilon)}c_k\ri\}.
\ena
It follows that
$$r_k^2=\f{\int_\Sigma he^{\langle\sigma_k,\zeta\rangle}dv_g}{\rho_k h(p)}e^{-c_k}\leq C
\exp\le\{-\f{1-3\epsilon}{2(1-\epsilon)}c_k\ri\}.$$
For any $0<\gamma<1/2$, we can take $\epsilon>0$ sufficiently small such that ${(1-3\epsilon)}/(2-2\epsilon)>\gamma$. As a consequence, we obtain
$r_k^2e^{\gamma c_k}\ra 0$ as $k\ra\infty$. $\hfill\Box$\\

Recall $\sigma_k=u_k\zeta$. In view of Lemma \ref{delta-L}, the Euler-Lagrange equation (\ref{E-Lk}) is transformed into
\be\label{E-L-uk}\Delta_gu_k=\rho_k\le(\f{he^{u_k}}{\mu_k}-\f{1}{|\Sigma|}\ri)-\lambda_{1,k}
\langle\tau_1,\zeta\rangle-(|\omega|^2+d^\ast\omega)u_k.\ee
For simplicity, in the isothermal coordinate systems
$\left(U, \phi_{k} ;\left\{y_{1}, y_{2}\right\}\right)$, we sometimes denote
$u \circ \phi_{k}^{-1}$ by $\widetilde{u}$ for a function $u$ in the sequel.
Define
$\psi_k(y)=c_k^{-1}\widetilde{u}_k(r_ky)$ and $\varphi_k(y)=\widetilde{u}_k(r_ky)-c_k$ for all
 $y\in \Omega_k=\{y\in\mathbb{R}^2:r_ky\in \mathbb{B}_{\delta }\}$. For the
convergence of $\psi_k$ and $\varphi_k$, one has
\begin{lemma}\label{bubble-1}
There hold up to a subsequence, $\psi_k\ra 1$ and $\varphi_k\ra
\varphi=-2\log(1+|y|^2/8)$ in $C^1_{\rm loc}(\mathbb{R}^2)$ as $k\ra\infty$.
\end{lemma}

\proof By (\ref{scale}), (\ref{E-L-uk}) and the definition of $g_k$, we calculate
\bea\nonumber
\Delta_{g_k}\psi_k(y)&=&\f{h(\phi_{k}^{-1}(r_ky))}{h(p)c_k}e^{\varphi_k(y)}-\f{r_k^2\rho_k}{|\Sigma|c_k}-c_k^{-1}r_k^2
\lambda_{1,k}\langle\tau_1, \zeta\rangle(\phi_k^{-1}(r_ky))\\\label{d-psi}
&&-c_k^{-1}r_k^2(|\omega|^2+d^\ast\omega)u_k
(\phi_k^{-1}(r_ky)),
\eea
where
$\Delta_{g_k}\psi_k(y)=-\exp\left(-{f_k}(r_ky)\right)\Delta_{\mathbb{R}^2}\psi_k(y).$
Since for any $R>0$,
$$\int_{\mathbb{B}_R} \le|r_k^2(|\omega|^2+d^\ast\omega)u_k(\phi_k^{-1}(r_ky))\ri|^2dy\leq Cr_k^2\int_\Sigma |\sigma_k|^2dv_g\leq
Cr_k^2c_k,$$
we have that $\Delta_{\mathbb{R}^2}\psi_k(y)$ converges to $0$ locally uniformly in $\mathbb{R}^2$. Note also that $\psi_k(y)\leq 1$ for all
$y\in\Omega_k$. Applying elliptic estimates to (\ref{d-psi}), we conclude that $\psi_k\ra \psi$ in
$C^1_{\rm loc}(\mathbb{R}^2)$, where $\psi$ is a harmonic function with $\psi(y)\leq \psi(0)=1$.
Then the Liouville theorem implies $\psi\equiv 1$ in $\mathbb{R}^2$.

Also we calculate
\bna\nonumber
\Delta_{g_k}\varphi_k(y)&=&\f{h(\phi_{k}^{-1}(r_ky))}{h(p)}e^{\varphi_k(y)}-\f{r_k^2\rho_k}{|\Sigma|}-r_k^2\lambda_{1,k}
\langle\tau_1,\zeta\rangle(\phi_k^{-1}(r_ky))\\\label{d-varphi}
&&-r_k^2(|\omega|^2+d^\ast\omega)u_k(\phi_k^{-1}(r_ky)).
\ena
Similarly as above, we use elliptic estimates to conclude that $\varphi_k\ra\varphi$ in $C^1_{\rm loc}(\mathbb{R}^2)$, where $\varphi$ is a solution of
$$\le\{\begin{array}{lll}
-\Delta_{\mathbb{R}^2}\varphi=e^{\varphi}\quad{\rm in}\quad \mathbb{R}^2\\[1.5ex]
\int_{\mathbb{R}^2}e^{\varphi(y)}dy<\infty.
\end{array}\ri.$$
It then follows from a classification theorem of Chen-Li \cite{Chen-Li} that
$\varphi(y)=-2\log(1+|y|^2/8)$ for all $y\in\mathbb{R}^2$, and thus
\be\label{integration}\int_{\mathbb{R}^2}e^{\varphi(y)}dy=8\pi,\ee
as desired. $\hfill\Box$
\begin{lemma}\label{concentration}
 $\mu_k^{-1}he^{u_k}dv_g$ converges to the dirac measure $\delta_p$ in the sense of measure.
\end{lemma}
\proof In view of (\ref{scale}) and Lemma \ref{bubble-1}, it is easy to get that
\bna\rho_k\int_{\phi_k^{-1}(\mathbb{B}_{Rr_k})}\mu_k^{-1}he^{u_k}dv_g&=&(1+o_k(1))\int_{\mathbb{B}_R}e^{\varphi_k(y)}dy\\
&=&(1+o_k(1))\int_{\mathbb{B}_R}e^{\varphi(y)}dy.\ena
This together with (\ref{integration}) gives
\be\label{d-1}\lim_{R\ra\infty}\lim_{k\ra\infty}\int_{\phi_k^{-1}(\mathbb{B}_{Rr_k})}\mu_k^{-1}he^{u_k}dv_g=1.\ee
Since $h(x)>0$ for all $x\in\Sigma$ and $\int_\Sigma \mu_k^{-1}he^{u_k}dv_g=1$, we obtain
\be\label{d-2}\lim_{R\ra\infty}\lim_{k\ra\infty}\int_{\Sigma\setminus \phi_k^{-1}(\mathbb{B}_{Rr_k})}\mu_k^{-1}he^{u_k}dv_g=0.\ee
Note also $x_k\ra p$ as $k\ra\infty$. It follows from (\ref{d-1}) and (\ref{d-2}) that
$$\lim_{k\ra\infty}\int_\Sigma \mu_k^{-1}\eta h e^{u_k}dv_g=\eta(p),\quad\forall\eta\in C^0(\Sigma).$$
This ends the proof of the lemma. $\hfill\Box$

\begin{lemma}\label{Green-0}
$\sigma_k$ converges to a Green section $\mathrm{G}$ weakly in $W^{1,q}(\Sigma,\mathcal{L})$, strongly in
$L^r(\Sigma,\mathcal{L})$ for all $r<\f{2q}{2-q}$, $1<q<2$, and in $C^1_{\rm loc}(\Sigma\setminus\{p\},\mathcal{L})$,
where $\mathrm{G}$ satisfies
\be\label{Green-equation}\le\{\begin{array}{lll}\Delta_{\mathcal{L}}{\rm G}=8\pi\le(\delta_p-\f{1}{|\Sigma|}\ri)\zeta-
   \lambda_{1}\tau_1
\\[1.5ex]
\lambda_{1}=8\pi\le(\langle\zeta,\tau_1\rangle(p)-\f{1}{|\Sigma|}\int_\Sigma\langle\zeta,\tau_1\rangle dv_g\ri)
\end{array}\ri.\ee
in the distributional sense.
\end{lemma}
\proof Denote a section sequence
$$\textsf{f}_k=\rho_k\le(\f{ he^{\langle\sigma_k,\zeta\rangle}}{\mu_k}-
\f{1}{|\Sigma|}\ri)\zeta-\lambda_{1,k}\tau_1.$$
We first prove that
\bna\|\sigma_k\|_{L^2(\Sigma,\mathcal{L})}=\le(\int_\Sigma|\sigma_k|^2dv_g\ri)^{1/2}\leq C.\ena
Suppose on the contrary $\|\sigma_k\|_{L^2(\Sigma,\mathcal{L})}\ra\infty$ as $k\ra\infty$.
In view of (\ref{E-Lk}), we have
\bna\label{eqn}\Delta_{\mathcal{L}}\widetilde{\sigma}_k=\widetilde{\textsf{f}}_k,\ena
where $\widetilde{\sigma}_k=\sigma_k/\|\sigma_k\|_{L^2(\Sigma,\mathcal{L})}$ and
$\widetilde{\textsf{f}}_k=\textsf{f}_k/\|\sigma_k\|_{L^2(\Sigma,\mathcal{L})}$. Since $\|\widetilde{\sigma}_k\|_{L^2(\Sigma,\mathcal{L})}=1$
and $\|\widetilde{\textsf{f}}_k\|_{L^1(\Sigma,\mathcal{L})}=o_k(1)$, we conclude from Proposition 7.1 in \cite{LiLiuYang} that
$\widetilde{\sigma}_k$ is bounded in $W^{1,q}(\Sigma,\mathcal{L})$ for any $1<q<2$. We assume with no loss of generality $\widetilde{\sigma}_k$
converges to $\widetilde{\sigma}$ weakly in $W^{1,q}(\Sigma,\mathcal{L})$, and $\widetilde{\sigma}$ is a distributional solution of
\be\label{harmonic}\Delta_{\mathcal{L}}\widetilde{\sigma}=0.\ee
Moreover $\widetilde{\sigma}$ belongs to $\mathcal{H}_1$ since $\widetilde{\sigma}_k\in\mathcal{H}_1$. Testing (\ref{harmonic})
with $\widetilde{\sigma}$ and using the Poincar\'e inequality
(Lemma \ref{Poincare}), we conclude $\|\widetilde{\sigma}\|_{L^2(\Sigma,\mathcal{L})}=0$, which contradicts
$\|\widetilde{\sigma}\|_{L^2(\Sigma,\mathcal{L})}=1$. Hence $\sigma_k$ must be bounded in $L^2(\Sigma,\mathcal{L})$. Again applying
Proposition 7.1 in \cite{LiLiuYang} to $\Delta_{\mathcal{L}}{\sigma_k}=\textsf{f}_k$, we have that $\sigma_k$ is bounded in $W^{1,q}(\Sigma,\mathcal{L})$
for any $1<q<2$. Then the weak compactness of $W^{1,q}(\Sigma,\mathcal{L})$ and the compact embedding
$W^{1,q}(\Sigma,\mathcal{L})\hookrightarrow L^r(\Sigma,\mathcal{L})$, $\forall r<{2q}/{(2-q)}$, lead to
the existence of some section ${\rm G}$ such that $\sigma_k$ converges to ${\rm G}$ weakly in $W^{1,q}(\Sigma,\mathcal{L})$ and strongly
in  $L^r(\Sigma,\mathcal{L})$, $\forall r<{2q}/{(2-q)}$, $1<q<2$. By (\ref{E-Lk}) and Lemma \ref{concentration}, ${\rm G}$ is a distributional solution
of (\ref{Green-equation}). While the convergence in $C^1_{\rm loc}(\Sigma\setminus\{p\},\mathcal{L})$ follows from elliptic estimates
on (\ref{E-Lk}). $\hfill\Box$\\

From Lemmas \ref{delta-L} and \ref{Green-0}, we know that there exists a unique $G_p\in W^{1,q}(\Sigma)$, $\forall 1<q<2$, such that ${\rm G}=G_p\zeta$ and in the distributional sense
\be\label{g-p}\Delta_gG_p+(|\omega|^2+d^\ast\omega)G_p=8\pi\le(\delta_p-\f{1}{|\Sigma|}\ri)
-\lambda_1\langle\tau_1,\zeta\rangle.\ee Let $r=r(x)={\rm d}_g(x,p)$ be the geodesic distance between $x$ and $p$.
Since $\Delta_g(G_p+4\log r)$ belongs to $L^q(\Sigma)$, the elliptic estimate implies that $G_p+4\log r$
belongs to $C^1(\Sigma)$. This together with the Sobolev embedding theorem leads to a local decomposition
\be \label{Green-p}G_p(x)=-4\log r+A_p+\eta(x),\ee
where $A_p$ is a constant, $\eta\in C^1(\Sigma)$ satisfies $\eta(p)=0$. As a consequence, the Green section ${\rm G}$
can be written as
\bna\label{Green-represent}{\rm G}(x)=\le(-4\log r+A_p+\eta(x)\ri)\zeta(x).\ena

\subsection{Lower bound estimate}
We now turn our attention to the lower bound estimate of $\mathcal{J}_{8\pi}$.
Firstly we calculate for any fixed $R>0$,
\bea\nonumber
 \int_{\phi_{k}^{-1}(\mathbb{B}_{Rr_k})}|\mathcal{D}\sigma_k|^2dv_g&=&\int_{\phi_{k}^{-1}(\mathbb{B}_{Rr_k})}
 |du_k\otimes\zeta+u_k\omega\otimes\zeta|^2dv_g\\\nonumber
&=&\int_{\phi_{k}^{-1}(\mathbb{B}_{Rr_k})}|du_k|^2dv_g+2\int_{\phi_{k}^{-1}(\mathbb{B}_{Rr_k})}u_k\langle du_k,\omega\rangle dv_g
\\&&+\int_{\phi_{k}^{-1}(\mathbb{B}_{Rr_k})}|u_k|^2|\omega|^2 dv_g,\label{bub-domain}
\eea
where $\phi_k:U\ra\mathbb{R}^2$ is the sequence of isothermal coordinate systems constructed before. By Lemma \ref{bubble-1},
we have $\widetilde{u}_k(z)=c_k-2\log(1+|z|^2/(8r_k^2))$ for all $z\in \mathbb{B}_{Rr_k}$. Hence
\bea\nonumber
\int_{\phi_{k}^{-1}(\mathbb{B}_{Rr_k})}|du_k|^2dv_g
&=&\int_{\mathbb{B}_{R}}|d\varphi_k(y)|^2dy\\&=&\int_{\mathbb{B}_{R}}|d\varphi(y)|^2dy+o_k(1)\nonumber\\
&=&16\pi\log{(1+\frac{R^{2}}{8})}-16\pi+o_k(1)+o_R(1),\label{main-ernergy}
\eea
and by Lemma \ref{r-0},
\be\label{o-1}
\int_{\phi_{k}^{-1}(\mathbb{B}_{Rr_k})}|u_k|^2|\omega|^2dv_g\leq C(1+o_k(1))R^2r_k^2c_k^2=o_k(1).
\ee
As a consequence we have by the H\"older inequality
\be\label{holder}\int_{\phi_{k}^{-1}(\mathbb{B}_{Rr_k})}u_k\langle du_k,\omega\rangle dv_g=o_k(1).\ee
Inserting (\ref{main-ernergy})-(\ref{holder}) into (\ref{bub-domain}), we obtain
\be\label{bubdom}
\int_{\phi_{k}^{-1}(\mathbb{B}_{Rr_k})}|\mathcal{D}\sigma_k|^2dv_g
=16\pi\log{\le(1+\frac{R^{2}}{8}\ri)}-16\pi+o_k(1)+o_R(1).
\ee

Secondly it follows from Lemma \ref{Green-0} that
\bea\nonumber
\int_{\Sigma\backslash \phi_{k}^{-1}(\mathbb{B}_\delta)}|\mathcal{D}\sigma_k|^2dv_g
&=&\int_{\Sigma\backslash \phi_{k}^{-1}(\mathbb{B}_\delta)}|\mathcal{D}\mathrm{G}|^2dv_g+o_k(1)\\\nonumber
&=&\int_{\Sigma\backslash \phi_{k}^{-1}(\mathbb{B}_\delta)}|dG_p|^2dv_g
+2\int_{\Sigma\backslash \phi_{k}^{-1}(\mathbb{B}_\delta)}G_p\langle dG_p,\omega\rangle dv_g
\\&&+\int_{\Sigma\backslash \phi_{k}^{-1}(\mathbb{B}_\delta)}|G_p|^2|\omega|^2 dv_g\label{G-est}
\eea
In view of (\ref{dist}), we conclude $d_g(x_k,\p\phi_k^{-1}(\mathbb{B}_\delta))=(1+O(\delta))\delta$. Thus there exists some constant
$c_1>0$ such that $B_{(1-c_1\delta)\delta}(x_k)\subset\phi_k^{-1}(\mathbb{B}_\delta)\subset B_{(1+c_1\delta)\delta}(x_k)$, provided that
$k$ is chosen sufficiently large. By (\ref{Green-p}), we have
\be\label{G-1}\int_{\phi_{k}^{-1}(\mathbb{B}_\delta)}G_p\langle dG_p,\omega\rangle dv_g=o(1)\ee
and
\be\label{G-2}
\int_{\phi_{k}^{-1}(\mathbb{B}_\delta)}|G_p|^2|\omega|^2 dv_g=o(1),
\ee
where $o(1)\ra 0$ as $k\ra\infty$ first, and then $\delta\ra 0$. Inserting (\ref{G-1}) and (\ref{G-2}) into (\ref{G-est}), we have by
integration by parts and (\ref{g-p}) that
\bea
\int_{\Sigma\backslash \phi_{k}^{-1}(\mathbb{B}_\delta)}|\mathcal{D}\sigma_k|^2dv_g
&=&\int_{\Sigma\backslash B_{(1+c_1\delta)\delta}(x_k)}|dG_p|^2dv_g
+\int_{\Sigma}\langle dG_p^2,\omega\rangle dv_g+\int_{\Sigma}G_p^2|\omega|^2 dv_g+o(1)\nonumber\\\nonumber
&\geq&\int_{\Sigma\backslash B_{(1+c_1\delta)\delta}(x_k)}G_p\Delta_g G_p dv_g-\int_{\partial B_{(1+c_1\delta)\delta}(x_k)}G_p\frac{\partial G_p}{\partial \nu}ds_g
\\&&+\int_{\Sigma}\langle dG_p^2,\omega\rangle dv_g+\int_{\Sigma}G_p^2|\omega|^2 dv_g+o(1)\nonumber\\
&=&-32\pi\log\delta+8\pi A_{p}-\frac{8\pi}{|\Sigma|}\int_{\Sigma}G_pdv_g
+o(1).\label{G-e}
\eea

Thirdly we estimate the energy of $\sigma_k$ on neck domains
$\phi_k^{-1}(\mathbb{B}_{\delta})\setminus \phi_k^{-1}(\mathbb{B}_{Rr_k})$.
Since $u_k$ converges to $G_p$ weakly in $W^{1,q}(\Sigma)$, $1<q<2$, and strongly  in $L^2(\Sigma)$,
$\phi_k^{-1}(\mathbb{B}_\delta)\subset B_{2\delta}(p)$ for sufficiently large $k$,
$\int_{B_{2\delta}(p)}G_p^2dv_g=o_\delta(1)$, and $\int_{B_{2\delta}(p)}G_p\langle dG_p,\omega\rangle dv_g=o_\delta(1)$,
there must be
$$\int_{\phi_k^{-1}(\mathbb{B}_\delta)\setminus\phi_k^{-1}(\mathbb{B}_{Rr_k})}u_k^2|\omega|^2dv_g\leq
C\int_{B_{2\delta}(p)}u_k^2dv_g=o(1)$$
and
$$\int_{\phi_k^{-1}(\mathbb{B}_\delta)\setminus\phi_k^{-1}(\mathbb{B}_{Rr_k})}u_k\langle du_k,\omega\rangle dv_g=o(1).$$
Therefore
\bea\nonumber
\int_{\phi_k^{-1}(\mathbb{B}_\delta)\setminus\phi_k^{-1}(\mathbb{B}_{Rr_k})}|\mathcal{D}\sigma_k|^2dv_g&=&
\int_{\phi_k^{-1}(\mathbb{B}_\delta)\setminus\phi_k^{-1}(\mathbb{B}_{Rr_k})}|du_k|^2dv_g+o(1)\\\nonumber
&=&\int_{\mathbb{B}_\delta\setminus\mathbb{B}_{Rr_k}}|d\widetilde{u}_k|^2dy+o(1)\\
&\geq&\int_{\mathbb{B}_\delta\setminus\mathbb{B}_{Rr_k}}|d\widetilde{u}_k^\ast|^2dy+o(1),\label{neck-1}
\eea
where $\widetilde{u}^\ast_k={\rm max}\{{\rm min}\{\widetilde{u}_k,a_k\},b_k\}$,
$a_k=\inf_{\partial \mathbb{B}_{Rr_k}}\widetilde{u}_k$, and $b_k=\sup_{\partial\mathbb{B}_\delta}\widetilde{u}_k$.
Let us now modify an argument of Wang-Liu \cite{Wang} (see also Li-Li \cite{LiLi} for the Toda system). It is easy to see that the infimum
\be\label{inf}\inf_{w|_{\partial \mathbb{B}_{Rr_k}}=a_k,\, w|_{\partial \mathbb{B}_\delta}=b_k}
\int_{\mathbb{B}_\delta\setminus\mathbb{B}_{Rr_k}}|dw|^2dy\ee
is attained by a harmonic function
\bna\label{solution}
H(y)=\frac{a_k-b_k}{\log (Rr_k)-\log \delta} \log |y|+\frac{a_k
 \log \delta-b_k \log (Rr_k)}{-\log (Rr_k)+\log \delta},\quad y\in \mathbb{B}_\delta\setminus\mathbb{B}_{Rr_k}.
\ena
As a consequence
\be\label{har-1}\int_{\mathbb{B}_\delta\setminus\mathbb{B}_{Rr_k}}|dH|^2dy=\f{2\pi(a_k-b_k)^2}{\log\delta-\log(Rr_k)}.\ee
Note that $\widetilde{u}_k^\ast=a_k$ on $\p\mathbb{B}_{Rr_k}$ and $\widetilde{u}_k^\ast=b_k$ on $\p\mathbb{B}_{\delta}$, provided
that $k$ is chosen sufficiently large.
Combining (\ref{neck-1})-(\ref{har-1}), we get by (\ref{scale})
\bea
\int_{\phi_k^{-1}(\mathbb{B}_\delta)\setminus\phi_k^{-1}(\mathbb{B}_{Rr_k})}|du_k|^2dv_g&\geq&
\f{2\pi(a_k-b_k)^2}{\log\delta-\log(Rr_k)}\nonumber\\
&=&\f{4\pi(c_k+a_k-c_k-b_k)^2}{c_k-\log\mu_k+2\log\delta-2\log R+\log\rho_k+\log h(p)}.\label{in-1}
\eea
In view of Lemmas \ref{bubble-1}, \ref{Green-0}, and the local representation (\ref{Green-p}), we derive
\be\label{dk}a_k-c_k-b_k=-2\log\le(1+\f{R^2}{8}\ri)+4\log\delta-A_p+o_k(1).\ee
It follows from Lemma \ref{r-0} that $c_k-\log\mu_k\geq \gamma c_k-\log\rho_k-\log h(p)$ for any $0<\gamma<1/2$,
in particular $\gamma\leq (c_k-\log\mu_k)/c_k\leq 1$, provided that $k$ is sufficiently large. This together with (\ref{dk}) leads  to a Taylor expansion of
the final term in (\ref{in-1}), namely
\bea\nonumber
\f{2\pi(a_k-b_k)^2}{\log\delta-\log(Rr_k)}&=&\f{4\pi c_k^2}{c_k-\log\mu_k}\le(1+\f{2(a_k-c_k-b_k)}{c_k}+
\f{(a_k-c_k-b_k)^2}{c_k^2}\ri)\\\nonumber
&&\times\le(1-\f{2\log\delta-2\log R+\log\rho_k+\log h(p)}{c_k-\log\mu_k}+O\le(\f{1}{(c_k-\log\mu_k)^2}\ri)\ri)\\
&=&\f{4\pi c_k^2}{c_k-\log\mu_k}+\f{8\pi (a_k-c_k-b_k)c_k}{c_k-\log\mu_k}-\f{4\pi c_k^2}{(c_k-\log\mu_k)^2}
(2\log\delta-2\log R\nonumber\\&&+\log\rho_k+\log h(p))+o_k(1).
\eea
Also we calculate
\bea\nonumber
\f{4\pi c_k^2}{c_k-\log\mu_k}&=&4\pi(c_k-\log\mu_k)\le(1+\f{\log\mu_k}{c_k-\log\mu_k}\ri)^2\\
&=&4\pi(c_k-\log\mu_k)\le(1-\f{\log\mu_k}{c_k-\log\mu_k}\ri)^2+16\pi\log\mu_k.\label{c2}
\eea
Combining (\ref{in-1})-(\ref{c2}) and (\ref{neck-1}), we obtain
\bea\nonumber
\int_{\phi_k^{-1}(\mathbb{B}_\delta)\setminus\phi_k^{-1}(\mathbb{B}_{Rr_k})}|\mathcal{D}\sigma_k|^2dv_g&\geq&
4\pi(c_k-\log\mu_k)\le(1-\f{\log\mu_k}{c_k-\log\mu_k}\ri)^2+16\pi\log\mu_k\\
&&+\f{8\pi (a_k-c_k-b_k)c_k}{c_k-\log\mu_k}-\f{4\pi c_k^2}{(c_k-\log\mu_k)^2}
(2\log\delta-2\log R\nonumber\\&&+\log\rho_k+\log h(p))+o(1).\label{neck-ener}
\eea
In view of (\ref{bubdom}), (\ref{G-e}), (\ref{neck-ener}) and $\rho_k\leq 8\pi$, there are the following estimates
\bea\nonumber
\mathcal{J}_{\rho_k}(\sigma_k)&=&\f{1}{2}\int_\Sigma|\mathcal{D}\sigma_k|^2dv_g+\f{\rho_k}{|\Sigma|}\int_\Sigma u_kdv_g-\rho_k\log\mu_k\\
&\geq&2\pi(c_k-\log\mu_k)\le(1-\f{\log\mu_k}{c_k-\log\mu_k}\ri)^2+\f{8\pi}{|\Sigma|}\int_\Sigma G_pdv_g
+\f{4\pi (a_k-c_k-b_k)c_k}{c_k-\log\mu_k}\nonumber\\&&-\f{2\pi c_k^2}{(c_k-\log\mu_k)^2}
(2\log\delta-2\log R+\log\rho_k+\log h(p))\nonumber\\&&-16\pi\log\delta+4\pi A_{p}-\frac{4\pi}{|\Sigma|}\int_{\Sigma}G_pdv_g
+8\pi\log{\le(1+\frac{R^{2}}{8}\ri)}-8\pi+o(1),\label{jk}
\eea
where $o(1)\ra 0$ as $k\ra\infty$ first, $R\ra\infty$ next, and $\delta\ra 0$ finally. Noting (\ref{dk}) and $\gamma c_k\leq c_k-\log\mu_k\leq c_k$
for any $\gamma<1/2$ and sufficiently large $k$, we have up to a subsequence
\be\label{1/2}\f{\log\mu_k}{c_k-\log\mu_k}=1+o_k(1).\ee
For otherwise, the righthand side of (\ref{jk}) tends to $+\infty$ as $k\ra\infty$. This implies $\mathcal{J}_{\rho_k}(\sigma_k)\ra+\infty$
as $k\ra\infty$, contradicting (\ref{les-1}). Note that (\ref{1/2}) implies
\be\label{lim-0}\f{c_k}{c_k-\log\mu_k}=2+o_k(1).\ee
Discarding the first term of the righthand side of (\ref{jk}), and passing to the limit
$k\ra\infty$ first, then $R\ra\infty$, $\delta\ra 0$ finally,  we have by using (\ref{dk}), (\ref{lim-0}) and Lemma \ref{lim-8pi} that
\be\label{low-1}
\inf_{\sigma\in\mathcal{H}_1}\mathcal{J}_{8\pi}(\sigma)\geq\f{4\pi}{|\Sigma|}
\int_\Sigma G_pdv_g-4\pi A_p-8\pi\log\pi-8\pi
-8\pi\log h(p).
\ee

We shall prove that (\ref{low-1}) is in fact an equality. It suffices to construct a sequence of sections ${\rm Q}_k\in\mathcal{H}_1$
such that $\mathcal{J}_{8\pi}({\rm Q}_k)$ converges to the righthand side of (\ref{low-1}). To this end, we set
\be\label{qk}q_k(x)=\le\{\begin{array}{lll}
c-2\log\le(1+\f{k^2r^2}{8}\ri),&x\in B_{R/k}(p)\\[1.2ex]
G_p(x)-\varrho(x)\eta(x),&x\in B_{2R/k}(p)\setminus B_{R/k}(p)\\[1.2ex]
G_p(x),&x\in\Sigma\setminus B_{2R/k}(p),
\end{array}\ri.\ee
where $\eta$ is given as in (\ref{Green-p}), $\varrho\in C_0^1(B_{2R/k}(p))$ satisfies $\varrho(x)\equiv 1$ for all $x\in B_{R/k}(p)$,
 $|d\varrho(x)|\leq 4k/R$ for all $x\in B_{2R/k}(p)$,
 $$c=2\log(1+R^2/8)-4\log R+4\log k+A_p,$$
 and $R=R(k)$ satisfies $R\ra+\infty$ and $R(\log R)/k\ra 0$ as $k\ra\infty$. One easily checks that $q_k$ belongs to $W^{1,2}(\Sigma)$ and
 \be\label{Q-L2}\int_\Sigma q_k^2dv_g=\int_\Sigma G_p^2dv_g+o_k(1).\ee
  Define
 \be\label{Qk}{\rm Q}_k=q_k\zeta-\int_\Sigma \langle q_k\zeta,\tau_1\rangle dv_g \tau_1,\ee
 where $\{\tau_1\}$ is a normal basis of $\mathcal{H}_0$. 
 Obviously ${\rm Q}_k\in\mathcal{H}_1$. In view of (\ref{Q-L2}) and $\mathcal{D}\tau_1=0$, we have
 \bea\nonumber
 \int_\Sigma|\mathcal{D}{\rm Q}_k|^2dv_g&=&\int_\Sigma|dq_k\otimes\zeta+q_k\omega\otimes\zeta|^2dv_g\\
 &=&\int_\Sigma\le(|dq_k|^2+q_k^2(|\omega|^2+d^\ast\omega)\ri)dv_g\nonumber\\
 &=&\int_\Sigma|dq_k|^2dv_g+\int_\Sigma G_p^2(|\omega|^2+d^\ast\omega)dv_g+o_k(1).\label{derive-Q}
 \eea
 By the definition of $q_k$ (see (\ref{qk}) above), we get
\be\label{bub-d}\int_{B_{R/k}(p)}|dq_k|^2dv_g=16\pi\log\le(1+\f{R^2}{8}\ri)-16\pi+o_k(1).\ee
 Moreover, in view of (\ref{qk}),
we have by integration by parts
\bea\nonumber
\int_{\Sigma\setminus B_{R/k}(p)}|dq_k|^2dv_g&=&\int_{\Sigma\setminus B_{R/k}(p)}|dG_p|^2dv_g+\int_{B_{2R/k}(p)\setminus B_{R/k}(p)}
|d(\varrho\eta)|^2dv_g\\\nonumber
&&-2\int_{B_{2R/k}(p)\setminus B_{R/k}(p)}\langle dG_p,d(\varrho\eta)\rangle dv_g\\\nonumber
&=&-\int_{\p B_{R/k}(p)}G_p\f{\p G_p}{\p\nu}ds_g+\int_{\Sigma\setminus B_{R/k}(p)}G_p\Delta_g G_pdv_g\\\nonumber
&&+\int_{B_{2R/k}(p)\setminus B_{R/k}(p)}|d(\varrho\eta)|^2dv_g-2\int_{B_{2R/k}(p)\setminus B_{R/k}(p)}\varrho\eta\Delta_gG_pdv_g\\
&&+2\int_{\p B_{R/k}(p)}\varrho\eta\f{\p G_p}{\p\nu}ds_g.\label{dom-2}
\eea
By (\ref{Green-p}) and (\cite{LiLiuYang}, Lemma 7.2), we have
\be\label{boundary}-\int_{\p B_{R/k}(p)}G_p\f{\p G_p}{\p\nu}ds_g=-32\pi\log R+32\pi\log k+8\pi A_p+o_k(1).\ee
Noting that
\be\label{perpend}\int_\Sigma G_p\langle\tau_1,\zeta\rangle dv_g=\int_\Sigma \langle\tau_1,{\rm G}\rangle dv_g=0,\ee
we also have by using (\ref{g-p}) and (\ref{Green-p}),
\be\int_{\Sigma\setminus B_{R/k}(p)}G_p\Delta_gG_pdv_g=-\f{8\pi}{|\Sigma|}\int_\Sigma G_pdv_g
-\int_\Sigma (|\omega|^2+d^\ast\omega)G_p^2dv_g+o_k(1).\label{out-domain}\ee
Using again (\ref{g-p}), (\ref{Green-p}) and (\cite{LiLiuYang}, Lemma 7.2), we know that all the last three terms on the righthand side of
(\ref{dom-2}) are infinitesimals as $k\ra\infty$.
In view of these infinitesimals,  (\ref{dom-2}) and (\ref{boundary}), we have by combining (\ref{derive-Q}), (\ref{bub-d}) and (\ref{out-domain}) that
\be
\int_\Sigma|\mathcal{D}{\rm Q}_k|^2dv_g=32\pi\log k-16\pi\log 8-16\pi+8\pi A_p-\f{8\pi}{|\Sigma|}\int_\Sigma G_pdv_g+o_k(1).\label{energy}
\ee
Next it is obvious to see
\bna\label{Q-int}\f{8\pi}{|\Sigma|}\int_\Sigma\langle {\rm Q}_k,\zeta\rangle dv_g=\f{8\pi}{|\Sigma|}\int_\Sigma
G_pdv_g+o_k(1).\ena
We are now computing the integral $\int_\Sigma he^{\langle{\rm Q}_k,\zeta\rangle}dv_g$. Fixing some small $\delta>0$, we write
\bea\nonumber
\int_\Sigma he^{\langle{\rm Q}_k,\zeta\rangle}dv_g&=&h(p)\int_{B_{R/k}(p)}e^{\langle{\rm Q}_k,\zeta\rangle}dv_g+
\int_{B_{R/k}(p)}(h-h(p))e^{\langle{\rm Q}_k,\zeta\rangle}dv_g\\
&&+\int_{B_\delta(p)\setminus B_{R/k}(p)}he^{\langle{\rm Q}_k,\zeta\rangle}dv_g+\int_{\Sigma\setminus B_\delta(p)}
he^{\langle{\rm Q}_k,\zeta\rangle}dv_g.\label{he}
\eea
Recalling (\ref{perpend}), there holds
$\int_\Sigma\langle q_k\zeta,\tau_1\rangle dv_g=o_k(1)$.
This together with (\ref{qk}) and (\ref{Qk}) gives
\bna h(p)\int_{B_{R/k}(p)}e^{\langle{\rm Q}_k,\zeta\rangle}dv_g&=&h(p)\int_{B_{R/k}(p)}e^{q_k+o_k(1)}dv_g\\
&=&(1+o_k(1))8\pi h(p)
e^{-2\log 8+2\log k+A_p}.\ena
While three integrals $\int_{B_{R/k}(p)}(h-h(p))e^{\langle{\rm Q}_k,\zeta\rangle}dv_g$,
$\int_{B_\delta(p)\setminus B_{R/k}(p)}he^{\langle{\rm Q}_k,\zeta\rangle}dv_g$ and
$\int_{\Sigma\setminus B_\delta(p)}
he^{\langle{\rm Q}_k,\zeta\rangle}dv_g$ are all $o_k(1)k^2$. These estimates together with (\ref{he}) lead to
\be\label{log}\log\int_\Sigma he^{\langle{\rm Q}_k,\zeta\rangle}dv_g=-\log 8+\log(\pi h(p))+2\log k+A_p+o_k(1).\ee
Combining (\ref{energy}) and (\ref{log}), we obtain
\bna
\mathcal{J}_{8\pi}({\rm Q}_k)&=&
\f{1}{2}\int_\Sigma|\mathcal{D}{\rm Q}_k|^2dv_g+\f{8\pi}{|\Sigma|}\int_\Sigma
\langle{\rm Q}_k,\zeta\rangle dv_g-8\pi\log \int_\Sigma he^{\langle{\rm Q}_k,\zeta\rangle}dv_g\\
&=&-8\pi\log\pi-8\pi-8\pi\log h(p)-4\pi A_p+\f{4\pi}{|\Sigma|}\int_\Sigma G_pdv_g+o_k(1).
\ena
This immediately leads to
\be\label{leq}\inf_{\sigma\in\mathcal{H}_1}\mathcal{J}_{8\pi}(\sigma)\leq-8\pi\log\pi-8\pi-8\pi\log h(p)-4\pi A_p+\f{4\pi}{|\Sigma|}\int_\Sigma G_pdv_g,\ee
which together with (\ref{low-1}) implies (\ref{lowerb}).

\subsection{Completion of the proof of Theorem \ref{T2}}

It is just a summary of the previous subsections. Our purpose is to prove that if $\mathcal{J}_{8\pi}$ has no
minimizer on $\mathcal{H}_1$,
that is,
\bna\label{minim}\mathcal{J}_{8\pi}(\sigma_0)\not=\inf_{\sigma\in\mathcal{H}_1}\mathcal{J}_{8\pi}(\sigma),\quad{\rm for\,\,all}\quad
\sigma_0\in
W^{1,2}(\Sigma,\mathcal{L}),\ena
then $\inf_{\sigma\in\mathcal{H}_1}\mathcal{J}_{8\pi}(\sigma)$ has an exact value as in (\ref{lowerb}).
First of all, let $\rho_k=8\pi-1/k$ for any $k\in\mathbb{N}^\ast$. By Theorem \ref{T1}, there exists a sequence of sections
$\sigma_k=u_k\zeta\in W^{1,2}(\Sigma,\mathcal{L})$ such that
$$\mathcal{J}_{\rho_k}(\sigma_k)=\inf_{\sigma\in\mathcal{H}_1}\mathcal{J}_{\rho_k}(\sigma).$$
By Lemma \ref{lim-8pi}, $\mathcal{J}_{\rho_k}(\sigma_k)\ra \inf_{\sigma\in\mathcal{H}_1}\mathcal{J}_{8\pi}(\sigma)$ as $k\ra\infty$.
Let $c_k=\max_\Sigma u_k=u_k(x_k)$. Then Lemma \ref{c-bdd} leads to $c_k\ra+\infty$ as $k\ra\infty$.
In view of Lemma \ref{Green-0}, there exists a unique Green section ${\rm G}=G_p\zeta$ such that $\sigma_k$ converges to ${\rm G}$ weakly
in $W^{1,q}(\Sigma,\mathcal{L})$ for any $1<q<2$, strongly in $L^r(\Sigma,\mathcal{L})$ for $r<2q/(2-q)$, $1<q<2$, and in $C^1_{\rm loc}
(\Sigma\setminus\{p\},\mathcal{L})$. Locally $G_p$ is written as $G_p(x)=-4\log r+A_p+o(1)$, where $o(1)\ra 0$ as $r=d_g(x,p)\ra 0$.
Finally we derive a lower bound (\ref{low-1}) of $\mathcal{J}_{8\pi}$ on $\mathcal{H}_1$. Also we construct a sequence of sections
${\rm Q}_k\in\mathcal{H}_1$ indicating that the infimum of $\mathcal{J}_{8\pi}$ on $\mathcal{H}_1$ satisfies (\ref{leq}). Combining
(\ref{low-1}) and (\ref{leq}), we obtain (\ref{lowerb}). This completes the proof of the theorem. $\hfill\Box$

   \section{Further assumption on the frame $\zeta$}\label{Sec4}
   We now turn to Theorem \ref{Thm-2}. A key observation is as follows.

   \begin{lemma}\label{parallel}
   Suppose that $\mathcal{D}\zeta=0$. Then ${\rm dim}\mathcal{H}_0=1$, in particular $\mathcal{H}_0={\rm span}\{\zeta\}$.
   Moreover if we set $\sigma=u\zeta$, then $\Delta_{\mathcal{L}}\sigma=(\Delta_gu)\zeta$.
   \end{lemma}
   \proof Since $\zeta$ is a global unit frame on $W^{1,2}(\Sigma,\mathcal{L})$, for any $\tau\in\mathcal{H}_0$, there exists a unique function
   $v$ such that $\tau=v\zeta$. Clearly
   $$0=\mathcal{D}\tau=dv\otimes\zeta+v\mathcal{D}\zeta=dv\otimes\zeta.$$
   This leads to $dv=0$, and thus $v\equiv C$ for some constant $C$. Hence $\tau=C\zeta$ and
   $\mathcal{H}_0={\rm span}\{\zeta\}$.

   In view of Lemma \ref{delta-L}, for any  $\sigma=u\zeta$, we have
   $$\Delta_{\mathcal{L}}\sigma=\le(\Delta_gu+(|\omega|^2+d^\ast\omega)u\ri)\zeta,$$
   where $\omega$ satisfies $\mathcal{D}\zeta=\omega\otimes\zeta$. Obviously $\omega=0$ since $\mathcal{D}\zeta=0$.
   Hence $\Delta_{\mathcal{L}}\sigma=(\Delta_gu)\zeta$.
   $\hfill\Box$\\

   {\it Proof of Theorem \ref{Thm-2}}.
   Suppose $\mathcal{D}\zeta=0$. If $\sigma=u\zeta\in\mathcal{H}_1$, then we have $\langle\sigma,\zeta\rangle=u$ and
   $$\int_\Sigma|\mathcal{D}\sigma|^2dv_g=\int_\Sigma|du|^2dv_g.$$
   In view of the definitions of functionals $J_\rho$ and $\mathcal{J}_\rho$ (see (\ref{Jrho}) and (\ref{J-funct}) above),
   there holds
   \be\label{funct=}\mathcal{J}_\rho(\sigma)=J_\rho(u).\ee
   Define a function space
   $$H_1=\le\{u\in W^{1,2}(\Sigma):\int_\Sigma udv_g=0\ri\}.$$

   {\it $(i)$ $\rho<8\pi$.}

   By Theorem \ref{T1}, there exists some $\sigma_\rho\in\mathcal{H}_1$ such that
   $\mathcal{J}_\rho(\sigma_\rho)=\inf_{\sigma\in\mathcal{H}_1}\mathcal{J}_\rho(\sigma)$. A straightforward calculation shows
   $\sigma_\rho$ satisfies the Euler-Lagrange equation (\ref{E-L}).
   In view of Lemma \ref{parallel}, one can easily check that (\ref{E-L}) reduces to (\ref{mean-eqn}). \\

   {\it (ii) $\rho=8\pi$.}

   Note that
   $$\inf_{\sigma\in\mathcal{H}_1}\mathcal{J}_{8\pi}(\sigma)=\inf_{u\in H_1}J_{8\pi}(u).$$
   Under the {\it Ding-Jost-Li-Wang condition} (\cite{DJLW}, Theorem 1.2), there exists some $u_0\in H_1$ such that
   $$J_{8\pi}(u_0)=\inf_{u\in H_1}J_{8\pi}(u).$$
   Let $\sigma_0=u_0\zeta$. Then $\sigma_0\in\mathcal{H}_1$ satisfies
   $$\mathcal{J}_{8\pi}(\sigma_0)=\inf_{\sigma\in\mathcal{H}_1}\mathcal{J}_{8\pi}(\sigma).$$
   Clearly $\sigma_0$ satisfies the Euler-Lagrange equation (\ref{mean-eqn}).\\

   {\it (iii) $\rho\in(8k\pi,8(k+1)\pi)$, $\forall k\in\mathbb{N}^\ast$.}

   It was proved by Djadli \cite{Djadli} that for any $\rho\in(8k\pi,8(k+1)\pi)$, $k\in\mathbb{N}^\ast$,
   $J_\rho$ has a critical point $u_\rho\in W^{1,2}(\Sigma)$. Let $\sigma_\rho=u_\rho\zeta$. In view of (\ref{funct=}),
   we conclude that $\sigma_\rho$ is a critical point of $\mathcal{J}_{\rho}$ on the Sobolev space
   $W^{1,2}(\Sigma,\mathcal{L})$. Obviously $\sigma_\rho$ is a solution of the mean field equation (\ref{mean-eqn}).\\

   Combining the above three assertions, we finish the proof of the theorem. $\hfill\Box$\\

   \noindent{\bf Acknowledgements}\\    

   We are very grateful to Professor Linlin Sun for discussion on this topic, which 
   enables us to improve a previous result ${\rm dim}\mathcal{H}_0<\infty$ to ${\rm dim}\mathcal{H}_0\leq 1$. 

\section*{Competing interests}
On behalf of all authors, the corresponding author states that there
 is no conflict of interest.
\section*{Data availability}
Data sharing not applicable to this article as no datasets were generated or analyzed during the current study.


\begin{thebibliography}{00}

   \bibitem{AJY}
   W. Ao, A. Jevnikar, W. Yang, On the Boundary Behavior for the Blow-up Solutions of the sinh-Gordon Equation and Rank $N$ Toda Systems in Bounded Domains, Int. Math. Res. Not. IMRN 2020 (2020), 9386-9419.

   \bibitem{Aubin}
   T. Aubin, Sur la function exponentielle, C. R. Acad. Sci. Paris S\'er. A-B 270 (1970), A1514-A1516.

   \bibitem{Aubin1982} T. Aubin, Nonlinear analysis on manifolds, Springer, 1982.

   \bibitem{Malchiodi2015} L. Battaglia, A. Jevnikar, A. Malchiodi, D. Ruiz, A general existence result for the Toda system on compact surfaces, Adv. Math. 285 (2015), 937-979.

   \bibitem{JDE20} L. Battaglia,  R. L\'opez-Soriano, A double mean field equation related to a curvature prescription problem, J. Differential Equations 269 (2020),
   2705-2740.

   \bibitem{Caffarelli}
   A. Caffarelli, Y. S. Yang, Vortex condensation in the Chern-Simons Higgs model: An existence theorem, Comm.
   Math. Phys. 168 (1995), 321-336.

   \bibitem{Chang-Yang1}
   A. Chang, P. Yang, Prescribing Gaussian curvature on $S^2$, Acta Math. Sci. Ser. B 159 (1987), 215-259.
   \bibitem{Chang-Yang2}
   A. Chang, P. Yang, Conformal deformation of metrics on $S^2$, J. Differential Geom. 23 (1988), 259-296.

   \bibitem{Chen-Li}
   W. Chen, C. Li, Prescribing Gaussian curvatures on surfaces with conical singularities,
   J. Geom. Anal. 1 (1991), 359-372.

   \bibitem{changesign} F. de Marchis, R. L\'opez-Soriano, D. Ruiz,  Compactness, existence and multiplicity for the singular mean
   field problem with sign-changing potentials, J. Math. Pures Appl. 115 (2018), 237-267.

   \bibitem{Djadli}
   Z. Djadli, Existence result for the mean field problem on Riemann surfaces of all genuses,
    Commun. Contemp. Math. 10 (2008), 205-220.

   \bibitem{DM}
   Z. Djadli, A. Malchiodi, Existence of conformal metrics with constant \textit{Q}-curvature, Ann.
   of Math. 168 (2008), 813-858.

   \bibitem{DJLW}
   W. Ding, J. Jost, J. Li,  G. Wang, The differential equation $\Delta u = 8\pi-8\pi he^u$ on a compact Riemann
   Surface, Asian J. Math. 1 (1997), 230-248.



   \bibitem{DJLW98}
   W. Ding, J. Jost, J. Li, G. Wang,  An analysis of the two-vortex case in the Chern-Simons Higgs model,
   Calc. Var. Partial Differential Equations 7 (1998), 87-97.

   \bibitem{DJLW99B}
   W. Ding, J. Jost, J. Li, G. Wang,  Multiplicity results for the two-vortex Chern-Simons Higgs model on the two-sphere,
   Comment. Math. Helv. 74 (1999), 118-142.

   \bibitem{DJLW99}
   W. Ding, J. Jost, J. Li,  G. Wang, Existence results for mean field equations, Ann. Inst. H. Poincar\'e Anal. Non Lin\'eaire. 16 (1999),
   653-666.

   \bibitem{Fontana}
   L. Fontana, Sharp borderline Sobolev inequalities on compact Riemannian manifolds, Comment. Math. Helv. 68 (1993), 415-454.


   \bibitem{KW}
   J. Kazdan, F. Warner, Curvature functions for compact $2$-manifolds, Ann. of Math. 99 (1974), 14-47.

   \bibitem{LiLi} J. Li, Y. Li, Solutions for Toda systems on Riemann surfaces, Ann. Sc. Norm. Super. Pisa Cl. Sci. 4 (2005), 703-728.

   \bibitem{LiLiuYang} Y. Li, P. Liu, Y. Yang, Moser-Trudinger inequalities on vector bundles over a compact Riemannian manifold of dimension 2, Calc. Var. Partial Differential Equations 28 (2007), 59-83.
   \bibitem{Li-Sun-Yang}
   J. Li, L. Sun, Y. Yang, Boundary value problem for the mean field equation on a compact Riemann surface, (2022), arXiv:2201.01544.

\bibitem{Li-Zhu}
   J. Li, C. Zhu, The convergence of the mean field type flow at a critical case, Calc. Var. Partial Differential Equations
   58 (2019), 60-78.


   \bibitem{Malchiodi-DCDS}
   A. Malchiodi,
   Topological methods for an elliptic equation with exponential nonlinearities,
   Discrete Contin. Dyn. Syst. 21 (2008), 277-294.
\bibitem{MN}
   A. Malchiodi, C. Ndiaye, Some existence results for the Toda system on closed surfaces,
    Atti Accad. Naz. Lincei Rend. Lincei Mat. Appl. 18 (2007), 391-412.

   \bibitem{MR}
   A. Malchiodi, D. Ruiz, A variational analysis of the Toda system on compact surfaces, Comm. Pure Appl. Math. 66 (2013), 332-371.

   \bibitem{Moser71}
   J. Moser, A sharp form of an inequality by N. Trudinger, Indiana Univ. Math. J. 20 (1971), 1077-1091.

   \bibitem{Moser73}
   J. Moser,  On a nonlinear problem in differential geometry, Dynamical systems (Proc. Sympos., Univ. Bahia, Salvador, 1971), Academic Press, New York, 1973, pp. 273-280.

   \bibitem{NT}
   M. Nolasco, G. Tarantello, On a sharp Sobolev-type inequality on two-dimensional compact manifolds, Arch. Ration. Mech. Anal. 145 (1998), 161-195.

   \bibitem{S}
   M. Struwe, The existence of surfaces of constant mean curvature with free boundaries,  Acta Math. Sci. Ser. B 160 (1988), 19-64.

   \bibitem{Sun-Wang-Yang}
   L. Sun, Y. Wang, Y. Yang,
   Existence results for a generalized mean field equation on a closed Riemann
   surface, (2021), arXiv: 2101.03859.

   \bibitem{Sun-Zhu}
   L. Sun, J. Zhu, Global existence and convergence of a flow to Kazdan-Warner equation with non-negative prescribed function, Calc. Var. Partial Differential Equations 60 (2021), Paper No. 42, 26 pp.

   \bibitem{Sun-Z}
   L. Sun, J. Zhu, Existence of Kazdan-Warner equation with sign-changing prescribed function, (2021), arXiv:
   2012.12840.

   \bibitem{Tarantello}G. Tarantello, Multiple condensate solutions for the Chern-Simons-Higgs theory, J. Math. Phys. 37 (1996), 3769-3796.

   \bibitem{Wang}
   M. Wang, Q. Liu, The equation $\Delta u+\nabla \phi \cdot \nabla u=8 \pi c\left(1-h e^{u}\right)$ on a Riemann surface,
   J. Partial Differ. Equ. 25 (2012), 335-355.


   \bibitem{YangYisong}
   Y. S. Yang, Solitons in field theory and nonlinear analysis, Springer Monographs in Mathematics, Springer, 2001.

   \bibitem{Yang-06} Y. Yang, A weighted form of Moser-Trudinger inequality on Riemannian surface, Nonlinear Anal. 65 (2006), 647-659.

   \bibitem{Yang2020}
   Y. Yang, Nonexistence of extremals for an inequality of Adimurthi-Druet on a closed Riemann surface, Sci. China Math. 63 (2020), 1627-1644.

   \bibitem{Yang-Zhu-Proc}
   Y. Yang, X. Zhu, A remark on a result of Ding-Jost-Li-Wang, Proc. Amer. Math. Soc. 145 (2017), 3953-3959.

   \bibitem{Yang-Zhu-Sci}Y. Yang, X. Zhu, Existence of solutions to a class of Kazdan-Warner equations on compact Riemannian surface, Sci. China Math. 61 (2018),  1109-1128.

\end{thebibliography}
\end{document}